\documentclass[a4paper,12pt]{article}

\usepackage{latexsym}
\usepackage{amssymb}
\usepackage{theorem}
\usepackage{amsmath}
\usepackage{amscd}
\usepackage{graphicx}
\usepackage{xcolor}
\usepackage{url}

\newtheorem{theorem}{Theorem}[section]
\newtheorem{corollary}[theorem]{Corollary}
\newtheorem{lemma}[theorem]{Lemma}
\newtheorem{example}[theorem]{Example}
\newtheorem{proposition}[theorem]{Proposition}
\newtheorem{remark}[theorem]{Remark}
\newtheorem{definition}[theorem]{Definition}

\newcommand{\demo}{\par\noindent{\it Proof. \/}\ }
\newcommand{\enD}{\hfill $\Box$\vspace{3truemm} \par}
\newcommand{\R}{\mathbb{R}}

\newcommand{\bn}{\mbox{\boldmath $n$}}
\newcommand{\bt}{\mbox{\boldmath $t$}}
\newcommand{\bs}{\mbox{\boldmath $s$}}
\newcommand{\ba}{\mbox{\boldmath $a$}}
\newcommand{\bb}{\mbox{\boldmath $b$}}

\newcommand{\be}{\mbox{\boldmath $e$}}

\newcommand{\bg}{\mbox{\boldmath $g$}}
\newcommand{\bv}{\mbox{\boldmath $v$}}
\newcommand{\bx}{\mbox{\boldmath $x$}}

\newcommand{\bw}{\mbox{\boldmath $w$}}

\begin{document}

\title{Bertrand framed surfaces in the Euclidean 3-space and its applications}

\author{Nozomi Nakatsuyama and Masatomo Takahashi}

\date{\today}

\maketitle
\begin{abstract} 
A framed surface is a smooth surface in the Euclidean space with a moving frame. 
By using the moving frame, we can define Bertrand framed surfaces as the same idea as Bertrand framed curves. 
Then we find the caustics and involutes as Bertrand framed surfaces.
As applications, we can directly define the caustics and involutes of framed surfaces, and give conditions that the caustics and involutes are inverse operations of framed surfaces like as those of Legendre curves. 
Moreover, a framed surface is one of the Bertrand framed surfaces if and only if another caustic of the involute exists, under conditions.
Furthermore, we find a new such operation, the so-called tangential direction framed surfaces.
\end{abstract}

\renewcommand{\thefootnote}{\fnsymbol{footnote}}

\section{Introduction}
Bertrand and Mannheim curves are classical objects in differential geometry (\cite{Aminov, Banchoff-Lovett, Berger-Gostiaux, Bertrand, doCarmo, HCIP, Izumiya-Takeuchi1, Kuhnel, Liu-Wang, Struik}). 
A Bertrand (respectively, Mannheim) curve in the Euclidean $3$-space is a space curve whose principal normal line is the same as the principal normal (respectively, bi-normal) line of another curve. 
By definition, another curve is a parallel curve with respect to the direction of the principal normal vector. 
In \cite{Honda-Takahashi-2020}, they investigated the conditions of the Bertrand and Mannheim curves of non-degenerate curves and framed curves. 
Moreover, we investigated the other cases, that is, a space curve whose tangent (or, principal normal, bi-normal) line is the same as the tangent (or, principal normal, bi-normal) line of another curve, respectively. 
We say that a Bertrand type curve if there exists such another curve. 
We investigated the existence conditions of Bertrand type curves in all cases in \cite{Nakatsuyama-Takahashi-1}. 
Moreover, we also investigated curves with singular points. 
As smooth curves with singular points, it is useful to use the framed curves in the Euclidean space (cf. \cite{Honda-Takahashi-2016}). 
We investigated the existence conditions of the Bertrand framed curves (Bertrand types of framed curves) in all cases in \cite{Nakatsuyama-Takahashi-1}. 
As a consequence, the involutes and circular evolutes of framed curves (cf. \cite{Honda-Takahashi-Preprint}) appear as the Bertrand framed curves. 
\par
On the other hand, a framed surface is a surface in Euclidean 3-space with a moving frame (cf. \cite{Fukunaga-Takahashi-2019}). 
Framed surfaces may have singular points. 
By using the moving frames, the basic invariants and the curvatures of framed surfaces are introduced in \cite{Fukunaga-Takahashi-2019}. 
By using the moving frame, we define Bertrand framed surfaces as the same idea as Bertrand framed curves. 
In this paper, we give existence conditions of Bertrand framed surfaces in all cases in \S 3. 
As a consequence, we find the caustics and involutes as Bertrand framed surfaces (Theorems \ref{ns-Bertrand}, \ref{nt-Bertrand}, \ref{sn-Bertrand} and \ref{tn-Bertrand}). 
For properties of the differential geometry of caustics see for example \cite{Arnold1, Gray, Izumiya-Takahashi-2008, Porteous, Teramoto-2019, Teramoto-2023}. 
As applications, we can directly define the caustics and involutes of framed surfaces, and give conditions that the caustics and involutes are inverse operations of framed surfaces (Theorem \ref{caustic-involute}) like as those of Legendre curves  (cf. \cite{Fukunaga-Takahashi-2015, Fukunaga-Takahashi-2016}) in \S 4. 
Moreover, a framed surface is one of the Bertrand framed surfaces if and only if another caustic of the involute exists, under conditions (Theorem \ref{caustic_involute_to_SS}).
Furthermore, we find a new such operation, the so-called tangential direction framed surfaces in \S 5. 
Finally, we give concrete examples of caustics, involutes and tangential direction framed surfaces in \S 6.
\par
We shall assume throughout the whole paper that all maps and manifolds are $C^{\infty}$ unless the contrary is explicitly stated.

\bigskip
\noindent
{\bf Acknowledgement}. 
The second author was partially supported by JSPS KAKENHI Grant Number JP 24K06728.

\section{Preliminaries}

Let $\R^3$ be the $3$-dimensional Euclidean space equipped with the inner product $\ba \cdot \bb = a_1 b_1 + a_2 b_2 + a_3 b_3$, 
where $\ba = (a_1, a_2, a_3)$ and $\bb = (b_1, b_2, b_3) \in \R^3$. 
The norm of $\ba$ is given by $\vert \ba \vert = \sqrt{\ba \cdot \ba}$ and the vector product is given by 
$$
\ba \times \bb={\rm det}\left(
\begin{array}{ccc}
\be_1 & \be_2 & \be_3 \\
a_1 & a_2 & a_3 \\
b_1 & b_2 & b_3 
\end{array}
\right),
$$
where $\be_1, \be_2, \be_3$ are the canonical basis on $\R^3$. 
Let $U$ be a simply connected domain of $\R^2$ and $S^2$ be the unit sphere in $\R^3$, that is, $S^2=\{\ba \in \R^3| |\ba|=1\}$.
We denote a $3$-dimensional smooth manifold $\{(\ba,\bb) \in S^2 \times S^2| \ba \cdot \bb=0\}$ by $\Delta$.

\begin{definition}\label{framed.surface}{\rm
We say that $(\bx,\bn,\bs):U \to \R^3 \times \Delta$ is a {\it framed surface} if $\bx_u (u,v) \cdot \bn (u,v)=\bx_v(u,v) \cdot \bn(u,v)=0$ for all $(u,v) \in U$, where $\bx_u(u,v)=(\partial \bx/\partial u)(u,v)$ and $\bx_v(u,v)=(\partial \bx/\partial v)(u,v)$. 
We say that $\bx:U \to \R^3$ is a {\it framed base surface} if there exists $(\bn,\bs):U \to \Delta$ such that $(\bx,\bn,\bs)$ is a framed surface.
}
\end{definition}
By definition, the framed base surface is a frontal. 
The definition and properties of frontals see \cite{Arnold1,Arnold2}.
On the other hand, the frontal is a framed base surface at least locally. 
In this paper, we consider framed base surfaces as singular surfaces. 

We denote $\bt(u,v)=\bn(u,v) \times \bs(u,v)$. 
Then $\{\bn(u,v),\bs(u,v),\bt(u,v)\}$ is a moving frame along $\bx(u,v)$.
Thus, we have the following systems of differential equations:
$$
\begin{pmatrix}
\bx_u \\
\bx_v
\end{pmatrix}
=
\begin{pmatrix}
a_1 & b_1 \\
a_2 & b_2 
\end{pmatrix}
\begin{pmatrix}
\bs \\
\bt
\end{pmatrix},
$$
$$
\begin{pmatrix}
\bn_u \\
\bs_u \\
\bt_u
\end{pmatrix}
=
\begin{pmatrix}
0 & e_1 & f_1 \\
-e_1 & 0 & g_1 \\
-f_1 & -g_1 & 0
\end{pmatrix}
\begin{pmatrix}
\bn \\
\bs \\
\bt
\end{pmatrix}
, \ 
\begin{pmatrix}
\bn_v \\
\bs_v \\
\bt_v
\end{pmatrix}
=
\begin{pmatrix}
0 & e_2 & f_2 \\
-e_2 & 0 & g_2 \\
-f_2 & -g_2 & 0
\end{pmatrix}
\begin{pmatrix}
\bn \\
\bs \\
\bt
\end{pmatrix},
$$
where $a_i,b_i,e_i,f_i,g_i:U \to \R, i=1,2$ are smooth functions and we call the functions {\it basic invariants} of the framed surface. 
We denote the above matrices by $\mathcal{G}, \mathcal{F}_1, \mathcal{F}_2$, respectively. 
We also call the matrices $(\mathcal{G}, \mathcal{F}_1, \mathcal{F}_2)$  {\it basic invariants} of the framed surface $(\bx,\bn,\bs)$. 
Note that $(u,v)$ is a singular point of $\bx$ if and only if ${\rm det}\ \mathcal{G}(u,v)=0$.

Since the integrability conditions $\bx_{uv}=\bx_{vu}$ and $\mathcal{F}_{2,u}-\mathcal{F}_{1,v}
=\mathcal{F}_1\mathcal{F}_2-\mathcal{F}_2\mathcal{F}_1$, 
the basic invariants should be satisfied the following conditions:
\begin{align}\label{integrability.condition}
\begin{cases}
a_{1v}-b_1g_2 = a_{2u}-b_2g_1, \\
b_{1v}-a_2g_1 = b_{2u}-a_1g_2, \\
a_1 e_2 + b_1 f_2 = a_2 e_1 + b_2 f_1,
\end{cases}
\begin{cases}
e_{1v}-f_1g_2 = e_{2u}-f_2g_1, \\
f_{1v}-e_2g_1 = f_{2u}-e_1g_2, \\
g_{1v}-e_1f_2 = g_{2u}-e_2f_1. 
\end{cases}
\end{align}
We have fundamental theorems for framed surfaces, that is, 
the existence and uniqueness theorem for the basic invariants of framed surfaces.

\begin{definition}\label{congruent.framed.surface}{\rm
Let $(\bx,\bn,\bs), (\widetilde{\bx},\widetilde{\bn},\widetilde{\bs}):U \to \R^3 \times \Delta$ be framed surfaces.
We say that $(\bx,\bn,\bs)$ and $(\widetilde{\bx},\widetilde{\bn},\widetilde{\bs})$ are {\it congruent as framed surfaces} if there exist a  constant rotation $A \in SO(3)$ and a translation $\ba \in \R^3$ such that 
$$
\widetilde{\bx}(u,v)=A(\bx(u,v))+\ba, \widetilde{\bn}(u,v)=A(\bn(u,v)), \widetilde{\bs}(u,v)=A(\bs(u,v)),
$$ 
for all $(u,v) \in U$.}
\end{definition}

\begin{theorem}[The Existence Theorem for framed surfaces]\label{existence.framed.surface}
Let $U$ be a simply \\ connected domain in $\R^2$ and let $a_i,b_i,e_i,f_i,g_i:U \to \R, i=1,2$ be smooth functions with the integrability conditions $(\ref{integrability.condition})$.
Then there exists a framed surface $(\bx,\bn,\bs):U \to \R^3 \times \Delta$ whose associated basic invariants is $(\mathcal{G}, \mathcal{F}_1, \mathcal{F}_2)$.
\end{theorem}

\begin{theorem}[The Uniqueness Theorem for framed surfaces]\label{uniqueness.framed.surface}
Let $(\bx,\bn,\bs),\\ (\widetilde{\bx},\widetilde{\bn},\widetilde{\bs}):U \to \R^3 \times \Delta$ be framed surfaces with the basic invariants  $(\mathcal{G},\mathcal{F}_1,\mathcal{F}_2),\\ (\widetilde{\mathcal{G}},\widetilde{\mathcal{F}}_1,\widetilde{\mathcal{F}}_2)$, respectively.
Then $(\bx,\bn,\bs)$ and $(\widetilde{\bx},\widetilde{\bn},\widetilde{\bs})$ are congruent as framed surfaces if and only if the basic invariants $(\mathcal{G},\mathcal{F}_1,\mathcal{F}_2)$ and $(\widetilde{\mathcal{G}},\widetilde{\mathcal{F}}_1,\widetilde{\mathcal{F}}_2)$ coincides.
\end{theorem}

Let $(\bx,\bn,\bs):U \to \R^3 \times \Delta$ be a framed surface with the basic invariants $(\mathcal{G},\mathcal{F}_1,\mathcal{F}_2)$.
For the moving frame $\{\bn,\bs,\bt\}$ along $\bx$, there is an ability.
We consider rotations and reflections of the vectors $\bs,\bt$. 
We denote
\begin{eqnarray*}
\begin{pmatrix}
\bs^\theta (u,v)\\
\bt^\theta (u,v)
\end{pmatrix}
=
\begin{pmatrix}
\cos \theta(u,v) & -\sin \theta(u,v) \\
\sin \theta(u,v) & \cos \theta(u,v) 
\end{pmatrix}
\begin{pmatrix}
\bs (u,v) \\
\bt (u,v)
\end{pmatrix},
\end{eqnarray*}
where $\theta:U \to \R$ is a smooth function. 
Then $\bn \times \bs^\theta=\bt^\theta$ and $\{\bn,\bs^\theta,\bt^\theta\}$ is also a moving frame along $\bx$.
It follows that $(\bx,\bn,\bs^\theta)$ is a framed surface.
We call the frame $\{\bn,\bs^\theta,\bt^\theta\}$ a {\it rotation frame}  by $\theta$ of the framed surface $(\bx,\bn,\bs)$.
We denote by $(\mathcal{G}^\theta,\mathcal{F}_1^\theta,\mathcal{F}_2^\theta)$ the basic invariants of $(\bx,\bn,\bs^\theta)$.
Moreover, we consider a moving frame $\{\bn^r,\bs^r,\bt^r\}=\{-\bn,\bt,\bs\}$ along $\bx$ and call it a {\it reflection frame} of the framed surface $(\bx,\bn,\bs)$. 
We denote by $(\mathcal{G}^r,\mathcal{F}_1^r,\mathcal{F}_2^r)$ the basic invariants of $(\bx,\bn^r,\bs^r)$.
By a direct calculation, we have the following.
\begin{proposition}\label{invariant_frame_change}
Under the above notations, we have the relationships between the basic invariants 
$(\mathcal{G},\mathcal{F}_1,\mathcal{F}_2)$ and $(\mathcal{G}^\theta,\mathcal{F}_1^\theta,\mathcal{F}_2^\theta)$, $(\mathcal{G}^r,\mathcal{F}_1^r,\mathcal{F}_2^r)$, respectively.
\par
$(1)$ For any smooth function $\theta:U \to \R$, 
\begin{eqnarray*}
\mathcal{G}^{\theta}
=\mathcal{G}
\begin{pmatrix}
\cos \theta & \sin \theta \\
-\sin \theta & \cos \theta
\end{pmatrix}
=
\begin{pmatrix}
a_1\cos \theta-b_1 \sin \theta & a_1 \sin \theta+b_1 \cos \theta \\
a_2\cos \theta-b_2 \sin \theta & a_2 \sin \theta+b_2 \cos \theta
\end{pmatrix},
\end{eqnarray*}
\begin{eqnarray*}
\mathcal{F}_1^{\theta}
=
\begin{pmatrix}
0 & e_1 \cos \theta - f_1 \sin \theta &  e_1 \sin \theta + f_1 \cos
 \theta \\
-e_1 \cos \theta + f_1 \sin \theta & 0 & g_1-\theta_u \\
-e_1 \sin \theta - f_1 \cos \theta & -g_1 + \theta_u & 0 
\end{pmatrix},
\end{eqnarray*}
\begin{eqnarray*}
\mathcal{F}_2^{\theta}
=
\begin{pmatrix}
0 & e_2 \cos \theta - f_2 \sin \theta &  e_2 \sin \theta + f_2 \cos
 \theta \\
-e_2 \cos \theta + f_2 \sin \theta & 0 & g_2-\theta_v \\
-e_2 \sin \theta - f_2 \cos \theta & -g_2 + \theta_v & 0 
\end{pmatrix}.
\end{eqnarray*}
\par
$(2)$ \begin{eqnarray*}
\mathcal{G}^{r}
=\mathcal{G}
\begin{pmatrix}
0 & 1 \\
1 & 0
\end{pmatrix}
=
\begin{pmatrix}
b_1 & a_1 \\
b_2 & a_2 
\end{pmatrix},
\mathcal{F}_1^{r}
=
\begin{pmatrix}
0 & - f_1  &  -e_1\\
f_1 & 0 & -g_1 \\
e_1 & g_1 & 0 
\end{pmatrix},
\mathcal{F}_2^{r}
=
\begin{pmatrix}
0 & -f_2 &  -e_2 \\
f_2 & 0 & -g_2 \\
e_2 & g_2 & 0 
\end{pmatrix}.
\end{eqnarray*}

\end{proposition}
Especially, we have 
\begin{eqnarray*}
\begin{pmatrix}
e_i^{\theta} \\
f_i^{\theta}
\end{pmatrix}
=
\begin{pmatrix}
\cos\theta & -\sin\theta \\
\sin \theta & \cos\theta 
\end{pmatrix}
\begin{pmatrix}
e_i \\
f_i
\end{pmatrix}, \ i=1,2.
\end{eqnarray*}

\begin{definition}\label{curvature.framed.surface}{\rm 
We define a smooth mapping $C^F=(J^F,K^F,H^F):U \to \R^3$ by 
\begin{eqnarray*}
J^F =
\det \begin{pmatrix} 
a_1 & b_1 \\ 
a_2 & b_2
\end{pmatrix}, 
K^F =
\det \begin{pmatrix} 
e_1 & f_1 \\ 
e_2 & f_2
\end{pmatrix}, \\
H^F = -\frac{1}{2}\left\{ 
\det \begin{pmatrix} 
a_1 & f_1 \\ 
a_2 & f_2
\end{pmatrix}
-
\det \begin{pmatrix} 
b_1 & e_1 \\ 
b_2 & e_2
\end{pmatrix}
\right\}.
\end{eqnarray*}
We call $C^F = (J^F,K^F,H^F)$ a {\it curvature of the framed surface}.
}
\end{definition}

The curvature is useful to recognize that the framed base surface is a front or not.
\begin{proposition}\label{Legendre.immersion}
Let $(\bx,\bn,\bs) : U \to \R^3 \times \Delta$ be a framed surface and $p \in U$. 
\par
$(1)$ Suppose that ${\rm rank} (d\bx)=1$ at $p$. 
Then $(\bx,\bn): U \to \R^3 \times S^2$ is a Legendre immersion around $p$ if and only if $H^F(p) \neq 0$.
\par
$(2)$ Suppose that ${\rm rank} (d\bx)=0$ at $p$. 
Then $(\bx,\bn): U \to \R^3 \times S^2$ is a Legendre immersion around $p$ if and only if $K^F(p) \neq 0$.
\end{proposition}

\section{Bertrand framed surfaces}

Let $(\bx,\bn,\bs)$ and $(\overline{\bx},\overline{\bn},\overline{\bs}):U \to \R^3 \times \Delta$ be framed surfaces. 
\begin{definition}\label{Bertrand-type-framed-surface}{\rm
We say that $(\bx,\bn,\bs)$ and $(\overline{\bx},\overline{\bn},\overline{\bs})$ are {\it $(\bv,\overline{\bw})$-mates} if there exists a smooth function $\lambda:U \to \R$ with $\lambda \not\equiv 0$ such that $\overline{\bx}(u,v)=\bx(u,v)+\lambda(u,v)\bv(u,v)$ and $\bv(u,v)= \overline{\bw}(u,v)$ for all $(u,v) \in U$, where $\bv$ and $\bw$ are $\bn, \bs$ or $\bt$, respectively. 
\par
We also say that  $(\bx,\bn,\bs)$ is a {\it $(\bv,\overline{\bw})$-Bertrand framed surface} (or, {\it $(\bv,\overline{\bw})$-Bertrand-\\Mannheim framed surface}) if there exists another framed surface $(\overline{\bx},\overline{\bn},\overline{\bs})$ such that $(\bx,\bn,\bs)$ and $(\overline{\bx},\overline{\bn},\overline{\bs})$ are $(\bv,\overline{\bw})$-mates.
}
\end{definition}

We clarify the notation $\lambda \not\equiv 0$. 
Throughout this paper, $\lambda \not\equiv 0$ means that $\{(u,v) \in U | \lambda(u,v) \not=0\}$ is a dense subset of $U$.
It follows that $\bx$ and $\overline{\bx}$ are different surfaces. 
Note that if $\lambda$ is constant, then $\lambda \not\equiv 0$ means that $\lambda$ is a non-zero constant.

Let $(\bx,\bn,\bs):U \to \R^3 \times \Delta$ be a framed surface with basic invariants $(\mathcal{G},\mathcal{F}_1,\mathcal{F}_2)$. 
We give existence conditions of Bertrand framed surfaces and basic invariants in all cases.

\begin{lemma}\label{nn-mate-lambda-const}
If $(\bx,\bn,\bs)$ and $(\overline{\bx},\overline{\bn},\overline{\bs}):U \to \R^3 \times \Delta$ are $(\bn,\overline{\bn})$-mates, then $\lambda$ is non-zero constant.
\end{lemma}
\demo
By definition, we have $\overline{\bx}(u,v)=\bx(u,v)+\lambda(u,v)\bn(u,v)$ and $\bn(u,v)=\overline{\bn}(u,v)$.
By differentiating, we have 
\begin{align*}
\overline{\bx}_u(u,v) &=\overline{a}_1(u,v)\overline{\bs}(u,v)+\overline{b}_1(u,v)\overline{\bt}(u,v)\\
&=(a_1(u,v)+\lambda(u,v)e_1(u,v))\bs(u,v)\\
&\quad+(b_1(u,v)+\lambda(u,v)f_1(u,v))\bt(u,v)
+\lambda_u(u,v)\bn(u,v),\\
\overline{\bx}_v(u,v) &=\overline{a}_2(u,v)\overline{\bs}(u,v)+\overline{b}_2(u,v)\overline{\bt}(u,v)\\
&=(a_2(u,v)+\lambda(u,v)e_2(u,v))\bs(u,v)\\
&\quad+(b_2(u,v)+\lambda(u,v)f_2(u,v))\bt(u,v)
+\lambda_v(u,v)\bn(u,v).
\end{align*}
Since $\overline{\bx}_u(u,v) \cdot \overline{\bn}(u,v)=\overline{\bx}_v(u,v) \cdot \overline{\bn}(u,v)=0$, $\lambda_u(u,v)=\lambda_v(u,v)=0$ for all $(u,v) \in U$. 
Therefore $\lambda$ is a constant. 
By $\lambda \not\equiv 0$, $\lambda$ is a non-zero constant.
\enD

\begin{theorem}\label{nn-Bertrand}
$(\bx,\bn,\bs): U \to \R^3 \times \Delta$ is always an $(\bn,\overline{\bn})$-Bertrand framed surface.
\end{theorem}
\demo
If we consider $(\overline{\bx},\overline{\bn},\overline{\bs}):U \to \R^3 \times \Delta$ by $(\overline{\bx},\overline{\bn},\overline{\bs})=(\bx+\lambda \bn,\bn,\bs)$, where $\lambda$ is a non-zero constant, then $(\overline{\bx},\overline{\bn},\overline{\bs})$ is a framed surface and $\bn=\overline{\bn}$. 
Hence, $(\bx,\bn,\bs): U \to \R^3 \times \Delta$ is an $(\bn,\overline{\bn})$-Bertrand framed surface.
\enD

By a direct calculation, we have the following (cf. \cite{Fukunaga-Takahashi-2019}).
\begin{proposition}\label{nn-Bertrand-basic-invariants}
Suppose that $(\bx,\bn,\bs)$ and $(\overline{\bx},\overline{\bn},\overline{\bs}):U \to \R^3 \times \Delta$ are $(\bn,\overline{\bn})$-mates, where $(\overline{\bx},\overline{\bn},\overline{\bs})=(\bx+\lambda \bn,\bn,\bs)$ and $\lambda$ is a non-zero constant. 
Then the basic invariants of $(\overline{\bx},\overline{\bn},\overline{\bs})$ are given by 
$$
\overline{\mathcal{G}}=\mathcal{G}+\lambda 
\begin{pmatrix}
e_1 & f_1\\
e_2 & f_2
\end{pmatrix}, \ 
\overline{\mathcal{F}}_1=\mathcal{F}_1, \ 
\overline{\mathcal{F}}_2=\mathcal{F}_2.
$$
\end{proposition}
\begin{remark}{\rm 
$(1)$ If $(\bx,\bn,\bs)$ and $(\overline{\bx},\overline{\bn},\overline{\bs}):U \to \R^3 \times \Delta$ are $(\bn,\overline{\bn})$-mates, then $\overline{\bx}$ is a parallel surface of $\bx$ (cf. \cite{Fukunaga-Takahashi-2019}).
\par 
$(2)$ On the moving frame of $(\overline{\bx},\overline{\bn},\overline{\bs})$, we can also take a rotation frame $\{\bn,\bs^\theta,\bt^\theta\}$ instead of $\{\bn,\bs,\bt\}$.
}
\end{remark}

\begin{theorem}\label{ns-Bertrand}
$(\bx,\bn,\bs):U \to \R^3 \times \Delta$ is an $(\bn,\overline{\bs})$-Bertrand framed surface if and only if there exist smooth functions $\lambda, \theta:U \to \R$ with $\lambda \not\equiv 0$ such that 
\begin{align}\label{ns-condition}
\begin{pmatrix}
a_1(u,v)+\lambda(u,v)e_1(u,v) & b_1(u,v)+\lambda(u,v)f_1(u,v))\\
a_2(u,v)+\lambda(u,v)e_2(u,v) & b_2(u,v)+\lambda(u,v)f_2(u,v))
\end{pmatrix}
\begin{pmatrix}
\sin \theta(u,v)\\
\cos \theta(u,v)
\end{pmatrix}
=
\begin{pmatrix}
0\\
0
\end{pmatrix}
\end{align}
for all $(u,v) \in U$.
\end{theorem}
\demo
Suppose that $(\bx,\bn,\bs):U \to \R^3 \times \Delta$ is an $(\bn,\overline{\bs})$-Bertrand framed surface. 
Then there exists $\lambda:U \to \R$ with $\lambda \not\equiv 0$ such that $\overline{\bx}(u,v)=\bx(u,v)+\lambda(u,v)\bn(u,v)$ and $\bn(u,v)=\overline{\bs}(u,v)$. 
By the same calculation of the proof of Theorem \ref{nn-mate-lambda-const}, we have 
$\overline{a}_1(u,v)=\lambda_u(u,v)$ and $\overline{a}_2(u,v)=\lambda_v(u,v)$. 
Moreover, since $\bn(u,v)=\overline{\bs}(u,v)$, there exists $\theta:U \to \R$ such that 
$$
\begin{pmatrix}
\overline{\bt}(u,v) \\
\overline{\bn}(u,v)
\end{pmatrix}
=
\begin{pmatrix}
\cos \theta(u,v) & -\sin \theta(u,v)\\
\sin \theta(u,v) & \cos \theta(u,v)
\end{pmatrix}
\begin{pmatrix}
\bs(u,v) \\
\bt(u,v)
\end{pmatrix}.
$$
Then we have 
\begin{align*}
\overline{b}_i(u,v)\cos \theta(u,v)&=a_i(u,v)+\lambda(u,v) e_i(u,v), \\
-\overline{b}_i(u,v)\sin \theta(u,v)&=b_i(u,v)+\lambda(u,v) f_i(u,v)
\end{align*}
for $i=1,2$. 
Therefore, we have 
$$
(a_i(u,v)+\lambda(u,v)e_i(u,v))\sin \theta(u,v)+(b_i(u,v)+\lambda(u,v)f_i(u,v))\cos \theta(u,v)=0
$$
for $i=1,2$ and all $(u,v) \in U$.
\par 
Conversely, suppose that there exist smooth functions $\lambda, \theta:U \to \R$ with $\lambda \not\equiv 0$ such that condition \eqref{ns-condition} satisfies. 
If we consider $(\overline{\bx}, \overline{\bn}, \overline{\bs}):U \to \R^3 \times \Delta$ as $(\bx+\lambda \bn,\sin \theta \bs+\cos \theta \bt, \bn)$, then 
we can show that $(\overline{\bx}, \overline{\bn}, \overline{\bs})$ is a framed surface. 
By definition, $(\bx,\bn,\bs)$ and $(\overline{\bx}, \overline{\bn}, \overline{\bs})$ are $(\bn,\overline{\bs})$-mates. 
\enD

\begin{proposition}\label{ns-Bertrand-basic-invariants}
Suppose that $(\bx,\bn,\bs)$ and $(\overline{\bx},\overline{\bn},\overline{\bs}):U \to \R^3 \times \Delta$ are $(\bn,\overline{\bs})$-mates, where $(\overline{\bx},\overline{\bn},\overline{\bs})=(\bx+\lambda \bn, \sin \theta \bs+\cos \theta \bt,\bn)$ and $\lambda, \theta:U \to \R$ are smooth functions satisfying $\lambda \not\equiv 0$ and condition \eqref{ns-condition}. 
Then the basic invariants of $(\overline{\bx},\overline{\bn},\overline{\bs})$ are given by 
\begin{align*}
\begin{pmatrix}
\overline{a}_1 & \overline{b}_1\\
\overline{a}_2 & \overline{b}_2
\end{pmatrix}
&=
\begin{pmatrix}
\lambda_u & (a_1+\lambda e_1) \cos \theta-(b_1+\lambda f_1) \sin \theta\\
\lambda_v & (a_2+\lambda e_2) \cos \theta-(b_2+\lambda f_2) \sin \theta
\end{pmatrix},\\
\begin{pmatrix}
\overline{e}_1 & \overline{f}_1 & \overline{g}_1\\
\overline{e}_2 & \overline{f}_2 & \overline{g}_2
\end{pmatrix}
&=
\begin{pmatrix}
-e_1 \sin \theta-f_1\cos \theta & \theta_u-g_1 & e_1\cos \theta-f_1\sin \theta\\
-e_2 \sin \theta-f_2\cos \theta & \theta_v-g_2 & e_2\cos \theta-f_2\sin \theta
\end{pmatrix}.
\end{align*}
\end{proposition}
\demo
By the proof of Theorem \ref{ns-Bertrand}, we have $\overline{a}_1=\lambda_u, \overline{a}_2=\lambda_v$ and 
$$
\overline{b}_1=(a_1+\lambda e_1) \cos \theta-(b_1+\lambda f_1) \sin \theta, 
\overline{b}_2=(a_2+\lambda e_2) \cos \theta-(b_2+\lambda f_2) \sin \theta.
$$
By differentiating $\overline{\bn}=\sin \theta \bs+\cos \theta \bt$, we have 
$\overline{\bn}_u=(\theta_u-g_1) \overline{\bt}+(-e_1\sin \theta-f_1\cos \theta) \overline{\bs}$ and $\overline{\bn}_v=(\theta_v-g_2) \overline{\bt}+(-e_2\sin \theta-f_2\cos \theta) \overline{\bs}$. 
Therefore, we have $\overline{e}_1=-e_1 \sin \theta -f_1 \cos \theta, \overline{f}_1=\theta_u-g_1$, $\overline{e}_2=-e_2 \sin \theta -f_2 \cos \theta$ and $\overline{f}_2=\theta_v-g_2$.
Moreover, by differentiating $\overline{\bt}=\cos \theta \bs-\sin \theta \bt$, we have $\overline{\bt}_u=(g_1-\theta_u) \overline{\bn}+(-e_1 \cos \theta+f_1 \sin \theta) \overline{\bs}$ and $\overline{\bt}_v=(g_2-\theta_v) \overline{\bn}+(-e_2 \cos \theta+f_2 \sin \theta) \overline{\bs}$. 
Therefore, we have $\overline{g}_1=e_1 \cos \theta-f_1 \sin \theta$ and  $\overline{g}_2=e_2 \cos \theta-f_2 \sin \theta$. 
\enD

\begin{remark}{\rm 
If $(\bx,\bn,\bs)$ and $(\overline{\bx},\overline{\bn},\overline{\bs}):U \to \R^3 \times \Delta$ are $(\bn,\overline{\bs})$-mates, then $\overline{\bx}$ is a caustic (an evolute or a focal surface) of $\bx$ (cf. \cite{Takahashi-Teramoto}). 
By condition \eqref{ns-condition}, we have
\begin{align}\label{det-cautics}
{\rm det} 
\begin{pmatrix}
a_1+\lambda e_1 & b_1+\lambda f_1\\
a_2+\lambda e_2 & b_2+\lambda f_2
\end{pmatrix}
=0.
\end{align}
It follows that $\lambda$ must be a solution of the equation $K^F \lambda^2-H^F \lambda+J^F=0$. 
It is easy to see that the converse does not hold in general in the case of $d\bx$ has a corank $2$ singular point, that is, condition \eqref{ns-condition} does not follows from \eqref{det-cautics}.
}
\end{remark}

\begin{theorem}\label{nt-Bertrand}
$(\bx,\bn,\bs):U \to \R^3 \times \Delta$ is an $(\bn,\overline{\bt})$-Bertrand framed surface if and only if there exist smooth functions $\lambda, \widetilde\theta:U \to \R$ with $\lambda \not\equiv 0$ such that 
\begin{align}\label{nt-condition}
\begin{pmatrix}
a_1(u,v)+\lambda(u,v)e_1(u,v) & b_1(u,v)+\lambda(u,v)f_1(u,v)\\
a_2(u,v)+\lambda(u,v)e_2(u,v) & b_2(u,v)+\lambda(u,v)f_2(u,v)
\end{pmatrix}
\begin{pmatrix}
-\cos \widetilde{\theta}(u,v)\\
\sin \widetilde{\theta}(u,v)
\end{pmatrix}
=
\begin{pmatrix}
0\\
0
\end{pmatrix}
\end{align}
for all $(u,v) \in U$.
\end{theorem}
\demo
Suppose that $(\bx,\bn,\bs):U \to \R^3 \times \Delta$ is an $(\bn,\overline{\bt})$-Bertrand framed surface. 
Then there exists $\lambda:U \to \R$ with $\lambda \not\equiv 0$ such that $\overline{\bx}(u,v)=\bx(u,v)+\lambda(u,v)\bn(u,v)$ and $\bn(u,v)=\overline{\bt}(u,v)$. 
By the same calculation of the proof of Lemma \ref{nn-mate-lambda-const}, we have 
$\overline{b}_1(u,v)=\lambda_u(u,v)$ and $\overline{b}_2(u,v)=\lambda_v(u,v)$. 
Moreover, since $\bn(u,v)=\overline{\bt}(u,v)$, there exists $\widetilde\theta:U \to \R$ such that 
$$
\begin{pmatrix}
\overline{\bn}(u,v) \\
\overline{\bs}(u,v)
\end{pmatrix}
=
\begin{pmatrix}
\cos \widetilde\theta(u,v) & -\sin \widetilde\theta(u,v)\\
\sin \widetilde\theta(u,v) & \cos \widetilde\theta(u,v)
\end{pmatrix}
\begin{pmatrix}
\bs(u,v) \\
\bt(u,v)
\end{pmatrix}.
$$
Then we have 
\begin{align*}
\overline{a}_i(u,v)\cos \widetilde\theta(u,v)&=b_i(u,v)+\lambda(u,v) f_i(u,v), \\
\overline{a}_i(u,v)\sin \widetilde\theta(u,v)&=a_i(u,v)+\lambda(u,v) e_i(u,v)
\end{align*}
for $i=1,2$. 
Therefore, we have 
$$
-(a_i(u,v)+\lambda(u,v)e_i(u,v))\cos \widetilde\theta(u,v)+(b_i(u,v)+\lambda(u,v)f_i(u,v))\sin \widetilde\theta(u,v)=0
$$
for $i=1,2$ and all $(u,v) \in U$.
\par 
Conversely, suppose that there exist smooth functions $\lambda, \widetilde{\theta}:U \to \R$ with $\lambda \not\equiv 0$ such that condition \eqref{nt-condition} satisfies. 
If we consider $(\overline{\bx}, \overline{\bn}, \overline{\bs}):U \to \R^3 \times \Delta$ as $(\bx+\lambda \bn,\cos \widetilde\theta \bs-\sin \widetilde\theta \bt, \sin \widetilde\theta \bs+\cos \widetilde\theta \bt)$, then 
we can show that $(\overline{\bx}, \overline{\bn}, \overline{\bs})$ is a framed surface. 
By definition, $(\bx,\bn,\bs)$ and $(\overline{\bx}, \overline{\bn}, \overline{\bs})$ are $(\bn,\overline{\bt})$-mates. 
\enD

\begin{proposition}\label{nt-Bertrand-basic-invariants}
Suppose that $(\bx,\bn,\bs)$ and $(\overline{\bx},\overline{\bn},\overline{\bs}):U \to \R^3 \times \Delta$ are $(\bn,\overline{\bt})$-mates, where $(\overline{\bx},\overline{\bn},\overline{\bs})=(\bx+\lambda \bn, \cos \widetilde\theta \bs-\sin \widetilde\theta \bt, \sin \widetilde\theta \bs+\cos \widetilde\theta \bt)$ and $\lambda, \theta:U \to \R$ are smooth functions satisfying $\lambda \not\equiv 0$ and condition \eqref{nt-condition}. 
Then the basic invariants of $(\overline{\bx},\overline{\bn},\overline{\bs})$ are given by 
\begin{align*}
\begin{pmatrix}
\overline{a}_1 & \overline{b}_1\\
\overline{a}_2 & \overline{b}_2
\end{pmatrix}
&=
\begin{pmatrix}
(a_1+\lambda e_1) \sin \widetilde\theta+(b_1+\lambda f_1) \cos \widetilde\theta & \lambda_u\\
(a_2+\lambda e_2) \sin \widetilde\theta+(b_2+\lambda f_2) \cos \widetilde\theta & \lambda_v
\end{pmatrix},\\
\begin{pmatrix}
\overline{e}_1 & \overline{f}_1 & \overline{g}_1\\
\overline{e}_2 & \overline{f}_2 & \overline{g}_2
\end{pmatrix}
&=
\begin{pmatrix}
g_1-\widetilde\theta_u & -e_1 \cos \widetilde\theta+f_1 \sin \widetilde\theta & -e_1 \sin \widetilde\theta-f_1 \cos \widetilde\theta\\
g_2-\widetilde\theta_v & -e_2 \cos \widetilde\theta+f_2 \sin \widetilde\theta & -e_2 \sin \widetilde\theta-f_2 \cos \widetilde\theta
\end{pmatrix}.
\end{align*}
\end{proposition}
\demo
By the proof of Theorem \ref{nt-Bertrand}, we have $\overline{b}_1=\lambda_u, \overline{b}_2=\lambda_v$ and 
$$
\overline{a}_1=(a_1+\lambda e_1) \sin \widetilde\theta+(b_1+\lambda f_1) \cos \widetilde\theta, 
\overline{a}_2=(a_2+\lambda e_2) \sin \widetilde\theta+(b_2+\lambda f_2) \cos \widetilde\theta.
$$
By differentiating $\overline{\bn}=\cos \widetilde\theta \bs-\sin \widetilde\theta \bt$, we have 
$\overline{\bn}_u=(g_1-\widetilde\theta_u) \overline{\bs}+(-e_1 \cos \widetilde\theta+\\f_1 \sin \widetilde\theta) \overline{\bt}$ and $\overline{\bn}_v=(g_2-\widetilde\theta_v) \overline{\bs}+(-e_2 \cos \widetilde\theta+f_2 \sin \widetilde\theta) \overline{\bt}$. 
Therefore, we have \\$\overline{e}_1=g_1-\widetilde\theta_u, \overline{f}_1=-e_1 \cos \widetilde\theta+f_1 \sin \widetilde\theta$, $\overline{e}_2=g_2-\widetilde\theta_v$ and $\overline{f}_2=-e_2 \cos \widetilde\theta+f_2 \sin \widetilde\theta$.\\
Moreover, by differentiating $\overline{\bs}=\sin \widetilde\theta \bs+\cos \widetilde\theta \bt$, we have $\overline{\bs}_u=-(g_1-\widetilde\theta_u) \overline{\bn}+(-e_1 \sin \widetilde\theta-f_1 \cos \widetilde\theta) \overline{\bt}$ and $\overline{\bs}_v=-(g_2-\widetilde\theta_v) \overline{\bn}+(-e_2 \sin \widetilde\theta-f_2 \cos \widetilde\theta) \overline{\bt}$. 
Therefore, we have $\overline{g}_1=-e_1 \sin \widetilde\theta-f_1 \cos \widetilde\theta$ and  $\overline{g}_2=-e_2 \sin \widetilde\theta-f_2 \cos \widetilde\theta$. 
\enD

\begin{theorem}\label{nt-Bertrand-equivalent}
$(\bx,\bn,\bs):U \to \R^3 \times \Delta$ is an $(\bn,\overline{\bt})$-Bertrand framed surface if and only $(\bx,\bn,\bs):U \to \R^3 \times \Delta$ is an $(\bn,\overline{\bs})$-Bertrand framed surface.
\end{theorem}
\demo
Suppose that $(\bx,\bn,\bs):U \to \R^3 \times \Delta$ is an $(\bn,\overline{\bt})$-Bertrand framed surface. By theorem \ref{nt-Bertrand}, there exist smooth functions $\lambda, \widetilde\theta:U \to \R$ with $\lambda \not\equiv 0$ such that the condition (\ref{nt-condition}). If $\widetilde\theta=\theta+\pi/2$, then we have $\sin\widetilde\theta=\cos \theta$ and $\cos\widetilde\theta=-\sin\theta$. Thus, we have the condition (\ref{ns-condition}). By theorem \ref{ns-Bertrand}, $(\bx,\bn,\bs)$ is an $(\bn,\overline{\bs})$-Bertrand framed surface. \par
Conversely, suppose that $(\bx,\bn,\bs):U \to \R^3 \times \Delta$ is an $(\bn,\overline{\bs})$-Bertrand framed surface. By theorem \ref{ns-Bertrand}, there exist smooth functions $\lambda, \theta:U \to \R$ with $\lambda \not\equiv 0$ such that the condition (\ref{ns-condition}). If $\theta=\widetilde\theta-\pi/2$, then we have $\sin\theta=-\cos\widetilde\theta$ and $\cos\theta=\sin\widetilde\theta$. Thus, we have the condition (\ref{nt-condition}). By theorem \ref{nt-Bertrand}, $(\bx,\bn,\bs)$ is an $(\bn,\overline{\bt})$-Bertrand framed surface.
\enD

\begin{theorem}\label{sn-Bertrand}
$(\bx,\bn,\bs):U \to \R^3 \times \Delta$ is an $(\bs,\overline{\bn})$-Bertrand framed surface if and only if ${\rm det}(\bb(u,v),\bg(u,v))=0$ for all $(u,v) \in U$ and $\lambda:U\to\R$ is given by  
$$
\lambda(u,v)=-\left(\int^u_{u_0} a_1(u,v) du+\int^v_{v_0} a_2(u_0,v) dv\right) +c
$$
for a point $(u_0,v_0) \in U$ and constant $c\in \R$ with $\lambda \not\equiv 0$.
\end{theorem}
\demo
Suppose that $(\bx,\bn,\bs):U \to \R^3 \times \Delta$ is an $(\bs,\overline{\bn})$-Bertrand framed surface. 
Then there exists $\lambda:U \to \R$ with $\lambda \not\equiv 0$ such that $\overline{\bx}(u,v)=\bx(u,v)+\lambda(u,v)\bs(u,v)$ and $\bs(u,v)=\overline{\bn}(u,v)$ for all $(u,v) \in U$. 
By differentiating, we have
\begin{align*}
&\overline{\bx}_u(u,v) =\overline{a}_1(u,v)\overline{\bs}(u,v)+\overline{b}_1(u,v)\overline{\bt}(u,v)\\
&=(a_1(u,v)+\lambda_u(u,v))\bs(u,v)\\
&\qquad+(b_1(u,v)+\lambda(u,v)g_1(u,v))\bt(u,v)-\lambda(u,v)e_1(u,v)\bn(u,v),\\
&\overline{\bx}_v(u,v) =\overline{a}_2(u,v)\overline{\bs}(u,v)+\overline{b}_2(u,v)\overline{\bt}(u,v)\\
&=(a_2(u,v)+\lambda_v(u,v))\bs(u,v)\\
&\qquad+(b_2(u,v)+\lambda(u,v)g_2(u,v))\bt(u,v)-\lambda(u,v)e_2(u,v)\bn(u,v).
\end{align*}
Since $\bs(u,v)=\overline{\bn}(u,v)$, we have $a_1(u,v)+\lambda_u(u,v)=0$ and $a_2(u,v)+\lambda_v(u,v)=0$ for all $(u,v) \in U$. 
It follows that $a_{1v}(u,v)=a_{2u}(u,v)$ and 
$$
\lambda(u,v)=-\left(\int^u_{u_0} a_1(u,v) du+\int^v_{v_0} a_2(u_0,v) dv \right) +c \not\equiv 0
$$
for a point $(u_0,v_0) \in U$ and constant $c\in \R$. 
By the integrability condition \eqref{integrability.condition}, we have ${\rm det}(\bb(u,v),\bg(u,v))=0$ for all $(u,v) \in U$. 
\par
Conversely, suppose that ${\rm det}(\bb(u,v),\bg(u,v))=0$, that is, $a_{1v}(u,v)=a_{2u}(u,v)$ for all $(u,v) \in U$. 
If we consider $(\overline{\bx},\overline{\bn},\overline{\bs}): U \to \R^3 \times \Delta$ as $(\bx+\lambda \bs,\bs,\bn)$, where 
$$
\lambda(u,v)=-\left(\int^u_{u_0} a_1(u,v) du+\int^v_{v_0} a_2(u_0,v) dv \right)+c,
$$
for a point $(u_0,v_0) \in U$ and constant $c\in \R$, 
then we can show that $(\overline{\bx},\overline{\bn},\overline{\bs})$ is a framed surface. 
By definition, $(\bx,\bn,\bs)$ and $(\overline{\bx}, \overline{\bn}, \overline{\bs})$ are $(\bs,\overline{\bn})$-mates. 
\enD

\begin{proposition}\label{sn-Bertrand-basic-invariants}
Suppose that $(\bx,\bn,\bs)$ and $(\overline{\bx},\overline{\bn},\overline{\bs}):U \to \R^3 \times \Delta$ are $(\bs,\overline{\bn})$-mates, where $(\overline{\bx},\overline{\bn},\overline{\bs})=(\bx+\lambda \bs, \bs,\bt)$ and $$
\lambda(u,v)=-\left(\int^u_{u_0} a_1(u,v) du+\int^v_{v_0} a_2(u_0,v) dv \right) +c \not\equiv 0
$$
for a point $(u_0,v_0) \in U$ and constant $c\in \R$.
Then the basic invariants of $(\overline{\bx},\overline{\bn},\overline{\bs})$ are given by 
\begin{align*}
\begin{pmatrix}
\overline{a}_1 & \overline{b}_1\\
\overline{a}_2 & \overline{b}_2
\end{pmatrix}
=
\begin{pmatrix}
b_1+\lambda g_1 & -\lambda e_1\\
b_1+\lambda g_2 & -\lambda e_2
\end{pmatrix}, \
\begin{pmatrix}
\overline{e}_1 & \overline{f}_1 & \overline{g}_1\\
\overline{e}_2 & \overline{f}_2 & \overline{g}_2
\end{pmatrix}
=
\begin{pmatrix}
g_1 & -e_1 & -f_1\\
g_2 & -e_2 & -f_2
\end{pmatrix}.
\end{align*}
\end{proposition}
\demo
By the proof of Theorem \ref{sn-Bertrand}, we have $\overline{a}_1=b_1+\lambda g_1, \overline{a}_2=b_1+\lambda g_2, \overline{b}_1=-\lambda e_1$ and $\overline{b}_2=-\lambda e_2$. 
By differentiating $\overline{\bn}=\bs$, we have 
$\overline{\bn}_u=-e_1\overline{\bt}+g_1\overline{\bs}$ and $\overline{\bn}_v=-e_2\overline{\bt}+g_2\overline{\bs}$. 
Therefore, we have $\overline{e}_1=g_1, \overline{f}_1=-e_1$, $\overline{e}_2=g_2$ and $\overline{f}_2=-e_2$.
Moreover, by differentiating $\overline{\bs}=\bt$, we have $\overline{\bs}_u=-f_1\overline{\bt}-g_1\overline{\bn}$ and $\overline{\bs}_v=-f_2\overline{\bt}-g_2\overline{\bn}$. 
Therefore, we have $\overline{g}_1=-f_1$ and  $\overline{g}_2=-f_2$. 
\enD

\begin{remark}{\rm 
If $(\bx,\bn,\bs)$ and $(\overline{\bx},\overline{\bn},\overline{\bs}):U \to \R^3 \times \Delta$ are $(\bs,\overline{\bn})$-mates, then we may consider $\overline{\bx}$ is one of involutes of $\bx$. 
}
\end{remark}

\begin{theorem}\label{ss-Bertrand}
$(\bx,\bn,\bs):U \to \R^3 \times \Delta$ is an $(\bs,\overline{\bs})$-Bertrand framed surface if and only if there exist smooth functions $\lambda, \theta:U \to \R$ with $\lambda \not\equiv 0$ such that 
\begin{align}\label{ss-condition}
\begin{pmatrix}
b_1(u,v)+\lambda(u,v)g_1(u,v) & \lambda(u,v)e_1(u,v) \\
b_2(u,v)+\lambda(u,v)g_2(u,v) & \lambda(u,v)e_2(u,v)
\end{pmatrix}
\begin{pmatrix}
\sin {\theta}(u,v)\\
-\cos {\theta}(u,v)
\end{pmatrix}
=
\begin{pmatrix}
0\\
0
\end{pmatrix}
\end{align}
\end{theorem}
for all $(u,v) \in U$.
\demo
Suppose that $(\bx,\bn,\bs):U \to \R^3 \times \Delta$ is an $(\bs,\overline{\bs})$-Bertrand framed surface. 
Then there exists $\lambda:U \to \R$ with $\lambda \not\equiv 0$ such that $\overline{\bx}(u,v)=\bx(u,v)+\lambda(u,v)\bs(u,v)$ and $\bs(u,v)=\overline{\bs}(u,v)$. 
By the same calculation of the proof of Theorem \ref{sn-Bertrand}, we have 
$\overline{a}_1=a_1+\lambda_u, \overline{a}_2=a_2+\lambda_v$. 
Moreover, since $\bs(u,v)=\overline{\bs}(u,v)$, there exists $\theta:U \to \R$ such that 
$$
\begin{pmatrix}
\overline{\bt}(u,v) \\
\overline{\bn}(u,v)
\end{pmatrix}
=
\begin{pmatrix}
\cos \theta(u,v) & -\sin \theta(u,v)\\
\sin \theta(u,v) & \cos \theta(u,v)
\end{pmatrix}
\begin{pmatrix}
\bt(u,v) \\
\bn(u,v)
\end{pmatrix}.
$$
Then we have 
$$
\overline{b}_i(u,v) \cos \theta(u,v)=b_i(u,v)+\lambda(u,v) g_i(u,v), \ \overline{b}_i(u,v) \sin \theta(u,v)=\lambda(u,v) e_i(u,v)
$$
for $i=1,2$. Therefore, we have 
$$
(b_i(u,v)+\lambda(u,v)g_i(u,v))\sin \theta(u,v)-\lambda(u,v)e_i(u,v)\cos \theta(u,v)=0
$$
for all $(u,v) \in U$ and $i=1,2$.
\par 
Conversely, suppose that there exist smooth functions $\lambda, \theta:U \to \R$ with $\lambda \not\equiv 0$ such that condition \eqref{ss-condition} satisfies. 
If we consider $(\overline{\bx}, \overline{\bn}, \overline{\bs}):U \to \R^3 \times \Delta$ as $(\bx+\lambda\bs, \sin \theta \bt+\cos \theta \bn, \bs)$, then 
we can show that $(\overline{\bx}, \overline{\bn}, \overline{\bs})$ is a framed surface. 
By definition, $(\bx,\bn,\bs)$ and $(\overline{\bx}, \overline{\bn}, \overline{\bs})$ are $(\bs,\overline{\bs})$-mates. 
\enD

\begin{proposition}\label{ss-Bertrand-basic-invariants}
Suppose that $(\bx,\bn,\bs)$ and $(\overline{\bx},\overline{\bn},\overline{\bs}):U \to \R^3 \times \Delta$ are $(\bs,\overline{\bs})$-mates, where $(\overline{\bx},\overline{\bn},\overline{\bs})=(\bx+\lambda \bs, \sin \theta \bt+\cos \theta \bn,\bs)$ and $\lambda, \theta:U \to \R$ are smooth functions satisfying $\lambda \not\equiv 0$ and condition \eqref{ss-condition}. 
Then the basic invariants of $(\overline{\bx},\overline{\bn},\overline{\bs})$ are given by 
\begin{align*}
\begin{pmatrix}
\overline{a}_1 & \overline{b}_1\\
\overline{a}_2 & \overline{b}_2
\end{pmatrix}
&=
\begin{pmatrix}
a_1+\lambda_u & (b_1+\lambda g_1) \cos \theta+\lambda e_1 \sin \theta\\
a_2+\lambda_v & (b_2+\lambda g_2) \cos \theta+\lambda e_2 \sin \theta
\end{pmatrix},\\
\begin{pmatrix}
\overline{e}_1 & \overline{f}_1 & \overline{g}_1\\
\overline{e}_2 & \overline{f}_2 & \overline{g}_2
\end{pmatrix}
&=
\begin{pmatrix}
-g_1 \sin \theta+e_1\cos \theta & \theta_u+f_1 & g_1\cos \theta+e_1\sin \theta\\
-g_2 \sin \theta+e_2\cos \theta & \theta_v+f_2 & g_2\cos \theta+e_2\sin \theta
\end{pmatrix}.
\end{align*}
\end{proposition}
\demo
By the proof of Theorem \ref{ss-Bertrand}, we have $\overline{a}_1=a_1+\lambda_u, \overline{a}_2=a_2+\lambda_v$ and 
$$
\overline{b}_1=(b_1+\lambda g_1) \cos \theta+\lambda e_1 \sin \theta, \ \overline{b}_2=(b_2+\lambda g_2) \cos \theta+\lambda e_2 \sin \theta.
$$
By differentiating $\overline{\bn}=\sin \theta \bt+\cos \theta \bn$, we have 
$\overline{\bn}_u=(\theta_u+f_1)\overline{\bt}+(-g_1 \sin \theta+e_1\cos \theta)\overline{\bs}$ and $\overline{\bn}_v=(\theta_v+f_2)\overline{\bt}+(-g_2 \sin \theta+e_2\cos \theta)\overline{\bs}$. 
Therefore, we have $\overline{e}_1=-g_1 \sin \theta+e_1\cos \theta, \overline{f}_1=\theta_u+f_1$, $\overline{e}_2=-g_2 \sin \theta+e_2\cos \theta$ and $\overline{f}_2=\theta_v+f_2$.
Moreover, by differentiating $\overline{\bt}=\cos\theta\bt-\sin\theta\bn$, we have $\overline{\bt}_u=-(\theta_u+f_1)\overline{\bn}-(g_1\cos \theta+e_1\sin \theta)\overline{\bs}$ and $\overline{\bt}_v=-(\theta_v+f_2)\overline{\bn}-(g_2\cos \theta+e_2\sin \theta)\overline{\bs}$. 
Therefore, we have $\overline{g}_1=g_1\cos \theta+e_1\sin \theta$ and  $\overline{g}_2=g_2\cos \theta+e_2\sin \theta$. 
\enD
\begin{remark}{\rm
If $(\bx,\bn,\bs)$ and $(\overline{\bx},\overline{\bn},\overline{\bs}):U \to \R^3 \times \Delta$ are $(\bs,\overline{\bs})$-mates, then 
$$
{\rm det} 
\begin{pmatrix}
b_1+\lambda g_1 & e_1\\
b_2+\lambda g_2 & e_2
\end{pmatrix}
=0
$$
by condition \eqref{ss-condition}. 
It follows that we have 
$$
{\rm det}(\bb (u,v),\be(u,v))+\lambda(u,v){\rm det}(\bg(u,v), \be(u,v))=0.
$$
If ${\rm det}(\be(u,v),\bg(u,v)) \not=0$, then $\lambda(u,v)={\rm det}(\bb(u,v),\be(u,v))/{\rm det}(\be(u,v),\bg(u,v))$. 
Hence, we have 
$$
\overline{\bx}(u,v)=\bx(u,v)+\frac{{\rm det}(\bb(u,v),\be(u,v))}{{\rm det}(\be(u,v),\bg(u,v))}\bs(u,v). 
$$
}
\end{remark}

\begin{theorem}\label{st-Bertrand}
$(\bx,\bn,\bs):U \to \R^3 \times \Delta$ is an $(\bs,\overline{\bt})$-Bertrand framed surface if and only if there exist smooth functions $\lambda, \widetilde\theta:U \to \R$ with $\lambda \not\equiv 0$ such that 
\begin{align}\label{st-condition}
\begin{pmatrix}
\lambda(u,v)e_1(u,v) & b_1(u,v)+\lambda(u,v)g_1(u,v)\\
\lambda(u,v)e_2(u,v) & b_2(u,v)+\lambda(u,v)g_2(u,v)
\end{pmatrix}
\begin{pmatrix}
\sin \widetilde{\theta}(u,v)\\
\cos \widetilde{\theta}(u,v)
\end{pmatrix}
=
\begin{pmatrix}
0\\
0
\end{pmatrix}
\end{align}
for all $(u,v) \in U$.
\end{theorem}
\demo
Suppose that $(\bx,\bn,\bs):U \to \R^3 \times \Delta$ is an $(\bs,\overline{\bt})$-Bertrand framed surface. 
Then there exists $\lambda:U \to \R$ with $\lambda \not\equiv 0$ such that $\overline{\bx}(u,v)=\bx(u,v)+\lambda(u,v)\bs(u,v)$ and $\bs(u,v)=\overline{\bt}(u,v)$. 
By the same calculation of the proof of Theorem \ref{ss-Bertrand}, we have 
$\overline{b}_1(u,v)=a_1(u, v)+\lambda_u(u,v)$ and $\overline{b}_2(u,v)=a_2(u, v)+\lambda_v(u,v)$. 
Moreover, since $\bs(u,v)=\overline{\bt}(u,v)$, there exists $\widetilde\theta:U \to \R$ such that 
$$
\begin{pmatrix}
\overline{\bn}(u,v) \\
\overline{\bs}(u,v)
\end{pmatrix}
=
\begin{pmatrix}
\cos \widetilde\theta(u,v) & -\sin \widetilde\theta(u,v)\\
\sin \widetilde\theta(u,v) & \cos \widetilde\theta(u,v)
\end{pmatrix}
\begin{pmatrix}
\bt(u,v) \\
\bn(u,v)
\end{pmatrix}.
$$
Then we have 
$$
\overline{a}_i(u,v)\cos \widetilde\theta(u,v)=-\lambda(u,v) e_i(u,v), \
\overline{a}_i(u,v)\sin \widetilde\theta(u,v)=b_i(u,v)+\lambda(u,v) g_i(u,v)
$$
for $i=1,2$. 
Therefore, we have 
$$
\lambda(u,v)e_i(u,v)\sin \widetilde\theta(u,v)+(b_i(u,v)+\lambda(u,v)g_i(u,v))\cos \widetilde\theta(u,v)=0
$$
for all $(u,v) \in U$ and $i=1,2$.
\par 
Conversely, suppose that there exist smooth functions $\lambda, \widetilde{\theta}:U \to \R$ with $\lambda \not\equiv 0$ such that condition \eqref{st-condition} satisfies. 
If we consider $(\overline{\bx}, \overline{\bn}, \overline{\bs}):U \to \R^3 \times \Delta$ as $(\bx+\lambda \bs,\cos \widetilde\theta \bt-\sin \widetilde\theta \bn, \sin \widetilde\theta \bt+\cos \widetilde\theta \bn)$, then 
we can show that $(\overline{\bx}, \overline{\bn}, \overline{\bs})$ is a framed surface. 
By definition, $(\bx,\bn,\bs)$ and $(\overline{\bx}, \overline{\bn}, \overline{\bs})$ are $(\bs,\overline{\bt})$-mates.
\enD
\begin{proposition}\label{st-Bertrand-basic-invariants}
Suppose that $(\bx,\bn,\bs)$ and $(\overline{\bx},\overline{\bn},\overline{\bs}):U \to \R^3 \times \Delta$ are $(\bs,\overline{\bt})$-mates, where $(\overline{\bx},\overline{\bn},\overline{\bs})=(\bx+\lambda \bs, \cos \widetilde\theta \bt-\sin \widetilde\theta \bn, \sin \widetilde\theta \bt+\cos \widetilde\theta \bn)$ and $\lambda, \widetilde\theta:U \to \R$ are smooth functions satisfying $\lambda \not\equiv 0$ and condition \eqref{st-condition}. 
Then the basic invariants of $(\overline{\bx},\overline{\bn},\overline{\bs})$ are given by 
\begin{align*}
\begin{pmatrix}
\overline{a}_1 & \overline{b}_1\\
\overline{a}_2 & \overline{b}_2
\end{pmatrix}
&=
\begin{pmatrix}
-\lambda e_1 \cos \widetilde\theta+(b_1+\lambda g_1) \sin \widetilde\theta & a_1+\lambda_u\\
-\lambda e_2 \cos \widetilde\theta+(b_2+\lambda g_2) \sin \widetilde\theta & a_2+\lambda_v 
\end{pmatrix},\\
\begin{pmatrix}
\overline{e}_1 & \overline{f}_1 & \overline{g}_1\\
\overline{e}_2 & \overline{f}_2 & \overline{g}_2
\end{pmatrix}
&=
\begin{pmatrix}
-\widetilde\theta_u-f_1 & -g_1\cos\widetilde\theta-e_1\sin\widetilde\theta & -g_1\sin\widetilde\theta+e_1\cos\widetilde\theta \\
-\widetilde\theta_v-f_2 & -g_2\cos\widetilde\theta-e_2\sin\widetilde\theta & -g_2\sin\widetilde\theta+e_2\cos\widetilde\theta
\end{pmatrix}.
\end{align*}
\end{proposition}
\demo
By the proof of Theorem \ref{st-Bertrand}, we have $\overline{b}_1=a_1+\lambda_u, \overline{b}_2=a_2+\lambda_v$, $\overline{a}_1=-\lambda e_1 \cos \widetilde\theta+(b_1+\lambda g_1) \sin \widetilde\theta$ and $ \overline{a}_2=-\lambda e_2 \cos \widetilde\theta+(b_2+\lambda g_2) \sin \widetilde\theta.
$
By differentiating $\overline{\bn}=\cos \widetilde\theta \bt-\sin \widetilde\theta \bn$, we have 
$\overline{\bn}_u=(-g_1\cos\widetilde\theta-e_1\sin\widetilde\theta)\overline{\bt}+(-\widetilde\theta_u-f_1)\overline{\bs}$ and $\overline{\bn}_v=(-g_2\cos\widetilde\theta-e_2\sin\widetilde\theta)\overline{\bt}+(-\widetilde\theta_v-f_2)\overline{\bs}$. 
Therefore, we have $\overline{e}_1=-\widetilde\theta_u-f_1, \overline{f}_1=-g_1\cos\widetilde\theta-e_1\sin\widetilde\theta$, $\overline{e}_2=-\widetilde\theta_v-f_2$ and $\overline{f}_2=-g_2\cos\widetilde\theta-e_2\sin\widetilde\theta$.
Moreover, by differentiating $\overline{\bs}=\sin\widetilde\theta\bt+\cos\widetilde\theta\bn$, we have $\overline{\bs}_u=(\widetilde\theta_u+f_1)\overline{\bn}+(-g_1\sin\widetilde\theta+e_1\cos\widetilde\theta)\overline{\bt}$ and $\overline{\bs}_v=(\widetilde\theta_v+f_2)\overline{\bn}+(-g_2\sin\widetilde\theta+e_2\cos\widetilde\theta)\overline{\bt}$. 
Therefore, we have $\overline{g}_1=-g_1\sin\widetilde\theta+e_1\cos\widetilde\theta$ and  $\overline{g}_2=-g_2\sin\widetilde\theta+e_2\cos\widetilde\theta$. 
\enD
\begin{theorem}\label{st-Bertrand-equivalent}
$(\bx,\bn,\bs):U \to \R^3 \times \Delta$ is an $(\bs,\overline{\bt})$-Bertrand framed surface if and only $(\bx,\bn,\bs):U \to \R^3 \times \Delta$ is an $(\bs,\overline{\bs})$-Bertrand framed surface.
\end{theorem}
\demo
Suppose that $(\bx,\bn,\bs):U \to \R^3 \times \Delta$ is an $(\bs,\overline{\bt})$-Bertrand framed surface. By Theorem \ref{st-Bertrand}, there exist smooth functions $\lambda, \widetilde\theta:U \to \R$ with $\lambda \not\equiv 0$ such that condition (\ref{st-condition}) satisfies. If $\widetilde\theta=\theta+\pi/2$, then we have $\sin\widetilde\theta=\cos \theta$ and $\cos\widetilde\theta=-\sin\theta$. Thus, we have condition (\ref{ss-condition}) satisfies. By Theorem \ref{ss-Bertrand}, $(\bx,\bn,\bs)$ is an $(\bs,\overline{\bs})$-Bertrand framed surface. \par
Conversely, suppose that $(\bx,\bn,\bs):U \to \R^3 \times \Delta$ is an $(\bs,\overline{\bs})$-Bertrand framed surface. By Theorem \ref{ss-Bertrand}, there exist smooth functions $\lambda, \theta:U \to \R$ with $\lambda \not\equiv 0$ such that condition (\ref{ss-condition}) satisfies. If $\theta=\widetilde\theta-\pi/2$, then we have $\sin\theta=-\cos\widetilde\theta$ and $\cos\theta=\sin\widetilde\theta$. Thus, we have condition (\ref{st-condition}) satisfies. By Theorem \ref{st-Bertrand}, $(\bx,\bn,\bs)$ is an $(\bs,\overline{\bt})$-Bertrand framed surface.
\enD

We can prove from Theorem \ref{tn-Bertrand} to Proposition \ref{tt-Bertrand-equivalent} by the similar calculations of proving of from Theorem \ref{sn-Bertrand} to Proposition \ref{st-Bertrand-equivalent}. 
Therefore, we omit the proof here. 

\begin{theorem}\label{tn-Bertrand}
$(\bx,\bn,\bs):U \to \R^3 \times \Delta$ is a $(\bt,\overline{\bn})$-Bertrand framed surface if and only if ${\rm det}(\bg(u,v),\ba(u,v))=0$ for all $(u,v) \in U$ and  and $\lambda:U\to\R$ is given by  
$$
\lambda(u,v)=-\left(\int^u_{u_0} b_1(u,v) du+\int^v_{v_0} b_2(u_0,v) dv\right) +c
$$
for a point $(u_0,v_0) \in U$ and constant $c\in \R$ with $\lambda\not\equiv0$.
\end{theorem}

\begin{proposition}\label{tn-Bertrand-basic-invariants}
Suppose that $(\bx,\bn,\bs)$ and $(\overline{\bx},\overline{\bn},\overline{\bs}):U \to \R^3 \times \Delta$ are $(\bt,\overline{\bn})$-mates, where $(\overline{\bx},\overline{\bn},\overline{\bs})=(\bx+\lambda \bt, \bt,\bn)$ and $$
\lambda(u,v)=-\left(\int^u_{u_0} b_1(u,v) du+\int^v_{v_0} b_2(u_0,v) dv \right) +c \not\equiv 0
$$
for a point $(u_0,v_0) \in U$ and constant $c\in \R$.
Then the basic invariants of $(\overline{\bx},\overline{\bn},\overline{\bs})$ are given by 
\begin{align*}
\begin{pmatrix}
\overline{a}_1 & \overline{b}_1\\
\overline{a}_2 & \overline{b}_2
\end{pmatrix}
=
\begin{pmatrix}
-\lambda f_1 & a_1-\lambda g_1\\
-\lambda f_2 & a_2-\lambda g_2
\end{pmatrix}, \
\begin{pmatrix}
\overline{e}_1 & \overline{f}_1 & \overline{g}_1\\
\overline{e}_2 & \overline{f}_2 & \overline{g}_2
\end{pmatrix}
=
\begin{pmatrix}
-f_1 & -g_1 & e_1\\
-f_2 & -g_2 & e_2
\end{pmatrix}.
\end{align*}
\end{proposition}

\begin{theorem}\label{ts-Bertrand}
$(\bx,\bn,\bs):U \to \R^3 \times \Delta$ is a $(\bt,\overline{\bs})$-Bertrand framed surface if and only if there exist smooth functions $\lambda, \theta:U \to \R$ with $\lambda \not\equiv 0$ such that 
\begin{align}\label{ts-condition}
\begin{pmatrix}
\lambda(u,v)f_1(u,v) & a_1(u,v)-\lambda(u,v)g_1(u,v)\\
\lambda(u,v)f_2(u,v) & a_2(u,v)-\lambda(u,v)g_2(u,v)
\end{pmatrix}
\begin{pmatrix}
-\sin {\theta}(u,v)\\
\cos {\theta}(u,v)
\end{pmatrix}
=
\begin{pmatrix}
0\\
0
\end{pmatrix}
\end{align}
for all $(u,v) \in U$.
\end{theorem}

\begin{proposition}\label{ts-Bertrand-basic-invariants}
Suppose that $(\bx,\bn,\bs)$ and $(\overline{\bx},\overline{\bn},\overline{\bs}):U \to \R^3 \times \Delta$ are $(\bt,\overline{\bs})$-mates, where $(\overline{\bx},\overline{\bn},\overline{\bs})=(\bx+\lambda \bt, \sin \theta \bn+\cos \theta \bs, \bt)$ and $\lambda, \theta:U \to \R$ are smooth functions satisfying $\lambda \not\equiv 0$ and condition \eqref{ts-condition}. 
Then the basic invariants of $(\overline{\bx},\overline{\bn},\overline{\bs})$ are given by 
\begin{align*}
\begin{pmatrix}
\overline{a}_1 & \overline{b}_1\\
\overline{a}_2 & \overline{b}_2
\end{pmatrix}
&=
\begin{pmatrix}
b_1+\lambda_u & -\lambda f_1 \cos \theta-(a_1-\lambda g_1) \sin \theta\\
b_2+\lambda_v & -\lambda f_2 \cos \theta-(a_2-\lambda g_2) \sin \theta
\end{pmatrix},\\
\begin{pmatrix}
\overline{e}_1 & \overline{f}_1 & \overline{g}_1\\
\overline{e}_2 & \overline{f}_2 & \overline{g}_2
\end{pmatrix}
&=
\begin{pmatrix}
f_1 \sin \theta+g_1\cos \theta & \theta_u-e_1 & -f_1\cos \theta+g_1\sin \theta\\
f_2 \sin \theta+g_2\cos \theta & \theta_v-e_2 & -f_2\cos \theta+g_2\sin \theta
\end{pmatrix}.
\end{align*}
\end{proposition}
\begin{theorem}\label{tt-Bertrand}
$(\bx,\bn,\bs):U \to \R^3 \times \Delta$ is a $(\bt,\overline{\bt})$-Bertrand framed surface if and only if there exist smooth functions $\lambda, \widetilde\theta:U \to \R$ with $\lambda \not\equiv 0$ such that 
\begin{align}\label{tt-condition}
\begin{pmatrix}
\lambda(u,v)f_1(u,v) & a_1(u,v)-\lambda(u,v)g_1(u,v)\\
\lambda(u,v)f_2(u,v) & a_2(u,v)-\lambda(u,v)g_2(u,v)
\end{pmatrix}
\begin{pmatrix}
\cos \widetilde{\theta}(u,v)\\
\sin \widetilde{\theta}(u,v)
\end{pmatrix}
=
\begin{pmatrix}
0\\
0
\end{pmatrix}
\end{align}
for all $(u,v) \in U$.
\end{theorem}

\begin{proposition}\label{tt-Bertrand-basic-invariants}
Suppose that $(\bx,\bn,\bs)$ and $(\overline{\bx},\overline{\bn},\overline{\bs}):U \to \R^3 \times \Delta$ are $(\bt,\overline{\bt})$-mates, where $(\overline{\bx},\overline{\bn},\overline{\bs})=(\bx+\lambda \bt, \cos \widetilde\theta \bn-\sin \widetilde\theta \bs, \sin \widetilde\theta \bn+\cos \widetilde\theta \bs)$ and $\lambda, \widetilde\theta:U \to \R$ are smooth functions satisfying $\lambda \not\equiv 0$ and condition \eqref{tt-condition}. 
Then the basic invariants of $(\overline{\bx},\overline{\bn},\overline{\bs})$ are given by 
\begin{align*}
\begin{pmatrix}
\overline{a}_1 & \overline{b}_1\\
\overline{a}_2 & \overline{b}_2
\end{pmatrix}
&=
\begin{pmatrix}
-\lambda f_1 \sin \widetilde\theta+(a_1-\lambda g_1) \cos \widetilde\theta & b_1+\lambda_u\\
-\lambda f_2 \sin \widetilde\theta+(a_2-\lambda g_2) \cos \widetilde\theta & b_2+\lambda_v 
\end{pmatrix},\\
\begin{pmatrix}
\overline{e}_1 & \overline{f}_1 & \overline{g}_1\\
\overline{e}_2 & \overline{f}_2 & \overline{g}_2
\end{pmatrix}
&=
\begin{pmatrix}
-\widetilde\theta_u-e_1 & f_1\cos\widetilde\theta-g_1\sin\widetilde\theta & f_1\sin\widetilde\theta+g_1\cos\widetilde\theta \\
-\widetilde\theta_v-e_2 & f_2\cos\widetilde\theta-g_2\sin\widetilde\theta & f_2\sin\widetilde\theta+g_2\cos\widetilde\theta
\end{pmatrix}.
\end{align*}
\end{proposition}

\begin{theorem}\label{tt-Bertrand-equivalent}
$(\bx,\bn,\bs):U \to \R^3 \times \Delta$ is a $(\bt,\overline{\bt})$-Bertrand framed surface if and only $(\bx,\bn,\bs):U \to \R^3 \times \Delta$ is a $(\bt,\overline{\bs})$-Bertrand framed surface.
\end{theorem}

\section{Caustics and involutes of framed surfaces}

The caustics (evolutes or focal surfaces) are classical object and it is well-known properties of caustics of regular surfaces (cf. \cite{Arnold1, Arnold2, Gray, Izumiya-book}).
Using Bertrand framed surfaces, we define caustics and involutes of framed surfaces directly. 
We denote 
$$\mathcal{F}(U,\R^3 \times \Delta):=\{(\bx,\bn,\bs)\in C^\infty(U,\R^3 \times \Delta)|(\bx,\bn,\bs) {\rm \ is\ a\ framed \ surface}\}.$$

Let $(\bx,\bn,\bs):U \to \R^3 \times \Delta$ be a framed surface with basic invariants $(\mathcal{G},\mathcal{F}_1,\mathcal{F}_2)$. 
\begin{definition}\label{caustic-framed-surface}{\rm 
$(1)$
The map $\mathcal{C}^{s}:\mathcal{F}(U,\R^3 \times \Delta)\to\mathcal{F}(U,\R^3 \times \Delta)$, $\mathcal{C}^{s}(\bx,\bn,\bs)=({\bx}^{n,s},{\bn}^{n,s},{\bs}^{n,s})$ is given by
\begin{align*}
&{\bx}^{n,s}(u,v)={\bx}(u,v)+{\lambda^{n,s}(u,v)}{\bn}(u,v),\\
&{\bn}^{n,s}(u,v)=\sin\theta^{n,s}(u,v)\bs(u,v)+\cos\theta^{n,s}(u,v)\bt(u,v),\\
&{\bs}^{n,s}(u,v)= \bn(u,v),
\end{align*}
where there exist smooth functions $\lambda^{n,s}, \theta^{n,s}:U \to \R$ such that 
{\small \begin{align}\label{ns-condition2}
\begin{pmatrix}
a_1(u,v)+\lambda^{n,s}(u,v) e_1(u,v) & b_1(u,v)+\lambda^{n,s}(u,v) f_1(u,v)\\
a_2(u,v)+\lambda^{n,s}(u,v) e_2(u,v) & b_2(u,v)+\lambda^{n,s}(u,v) f_2(u,v)
\end{pmatrix}
\begin{pmatrix}
\sin\theta^{n,s}(u,v) \\
\cos\theta^{n,s}(u,v)
\end{pmatrix}
=
\begin{pmatrix}
0 \\
0
\end{pmatrix}
\end{align}}
for all $(u,v) \in U$. 
Then we say that ${\bx}^{n,s}:U \to \R^3$ is a {\it caustic} of the framed surface $(\bx,\bn,\bs)$. \par
$(2)$
The map $\mathcal{C}^{t}:\mathcal{F}(U,\R^3 \times \Delta)\to\mathcal{F}(U,\R^3 \times \Delta)$, $\mathcal{C}^{t}(\bx,\bn,\bs)=({\bx}^{n,t},{\bn}^{n,t},{\bs}^{n,t})$ is given by
\begin{align*}
&{\bx}^{n,t}(u,v)={\bx}(u,v)+{\lambda^{n,t}(u,v)}{\bn}(u,v),\\
&{\bn}^{n,t}(u,v)=\cos \theta^{n,t}(u,v) \bs(u,v)-\sin \theta^{n,t}(u,v) \bt(u,v),\\
&{\bs}^{n,t}(u,v)=\sin \theta^{n,t}(u,v) \bs(u,v)+\cos \theta^{n,t}(u,v) \bt(u,v),
\end{align*}
where there exist smooth functions $\lambda^{n,t}, \theta^{n,t}:U \to \R$ such that 
{\small \begin{align}\label{nt-condition2}
\begin{pmatrix}
a_1(u,v)+\lambda^{n,t}(u,v)e_1(u,v) & b_1(u,v)+\lambda^{n,t}(u,v)f_1(u,v)\\
a_2(u,v)+\lambda^{n,t}(u,v)e_2(u,v) & b_2(u,v)+\lambda^{n,t}(u,v)f_2(u,v)
\end{pmatrix}
\begin{pmatrix}
-\cos\theta^{n,t}(u,v) \\
\sin\theta^{n,t}(u,v)
\end{pmatrix}
=
\begin{pmatrix}
0 \\
0
\end{pmatrix}
\end{align}}
for all $(u,v) \in U$. 
Then we say that ${\bx}^{n,t}:U \to \R^3$ is a {\it caustic} of the framed surface $(\bx,\bn,\bs)$. 
}
\end{definition}

\begin{remark}{\rm 
$(1)$ The caustic ${\bx}^{n,s}$ (respectively, ${\bx}^{n,t}$) is corresponding to the $(\bn, \overline{\bs})$ (respectively, $(\bn, \overline{\bt})$)-Bertrand framed surface. 
\par
$(2)$ By a direct calculation, we have ${\bt}^{n,s}(u,v)=\cos\theta^{n,s}(u,v)\bs(u,v)-\sin\theta^{n,s}(u,v)\bt(u,v)$ and ${\bt}^{n,t}(u,v)=\bn(u,v)$. 
\par
$(3)$ Suppose that there exist smooth functions $\lambda^{n,s}, \theta^{n,s}:U \to \R$ such that the condition (\ref{ns-condition2}) satisfies. If we take smooth functions $\lambda^{n,t}, \theta^{n,t}:U \to \R$ by $\lambda^{n,t}=\lambda^{n,s}$ and $\theta^{n,t}=\theta^{n,s}+\pi/2$, then the condition (\ref{nt-condition2}) is satisfied (cf. Theorem \ref{nt-Bertrand-equivalent}). The reflection frame of $\mathcal{C}^{s}(\bx,\bn,\bs)$ is corresponding to the moving frame of $\mathcal{C}^{t}(\bx,\bn,\bs)$. It follows that the map $\mathcal{C}^{t}$ is given by $\mathcal{C}^{t}(\bx,\bn,\bs)=\mathcal{C}^{s}(\bx,-\bn,\bt)$.
}
\end{remark}

\begin{definition}\label{involute-framed-surface}{\rm
$(1)$ Suppose that ${\rm det}(\bb(u,v),\bg(u,v))=0$ for all $(u,v) \in U$ and $(u_0,v_0) \in U$. 
The map $\mathcal{I}^{s}:\mathcal{F}(U,\R^3 \times \Delta)\to\mathcal{F}(U,\R^3 \times \Delta)$, $\mathcal{I}^{s}(\bx,\bn,\bs)=({\bx}^{s,n},{\bn}^{s,n},{\bs}^{s,n})$ is given by
\begin{align*}
&{\bx}^{s,n}(u,v)={\bx}(u,v)+\lambda^{s,n}(u,v)\bs(u,v),\\
&{\bn}^{s,n}(u,v)=\bs(u,v),\\
&{\bs}^{s,n}(u,v)=\cos \theta^{s,n}(u,v) \bt(u,v)-\sin \theta^{s,n}(u,v) \bn(u,v),
\end{align*}
where $\theta^{s,n}:U \to \R$ is a smooth function and $\lambda^{s,n}:U\to\R$ is given by 
$$
\lambda^{s,n}(u,v)=-\left(\int^u_{u_0} a_1(u,v) du+\int^v_{v_0} a_2(u_0,v) dv\right).
$$ 
Then we say that ${\bx}^{s,n}:U \to \R^3$ is an {\it involute with respect to $\bs$} at $(u_0,v_0) \in U$ of the framed surface $(\bx,\bn,\bs)$. \par
$(2)$ Suppose that ${\rm det}(\ba(u,v),\bg(u,v))=0$ for all $(u,v) \in U$ and $(u_0,v_0) \in U$. 
The map $\mathcal{I}^{t}:\mathcal{F}(U,\R^3 \times \Delta)\to\mathcal{F}(U,\R^3 \times \Delta)$, $\mathcal{I}^{t}(\bx,\bn,\bs)=({\bx}^{t,n},{\bn}^{t,n},{\bs}^{t,n})$ is given by
\begin{align*}
&{\bx}^{t,n}(u,v)={\bx}(u,v)+\lambda^{t,n}(u,v)\bt(u,v),\\
&{\bn}^{t,n}(u,v)=\bt(u,v),\\ 
&{\bs}^{t,n}(u,v)=\cos\theta^{t,n}(u,v) \bn(u,v)-\sin\theta^{t,n}(u,v) \bs(u,v),
\end{align*}
where $\theta^{t,n}:U \to \R$ is a smooth function and $\lambda^{t,n}:U\to\R$ is given by 
$$
\lambda^{t,n}(u,v)=-\left(\int^u_{u_0} b_1(u,v) du+\int^v_{v_0} b_2(u_0,v) dv\right).
$$ 
Then we say that ${\bx}^{t,n}:U \to \R^3$ is an {\it involute with respect to $\bt$} at $(u_0,v_0) \in U$ of the framed surface $(\bx,\bn,\bs)$. 
}
\end{definition}

\begin{remark}{\rm
$(1)$ The involute ${\bx}^{s,n}$ (respectively, ${\bx}^{t,n}$) is corresponding to the $(\bs, \overline{\bn})$ (respectively, $(\bt, \overline{\bn})$)-Bertrand framed surface under the condition $\theta^{s,n}=0$ (respectively, $\theta^{t,n}=0$). 
However, we consider a framed rotation of the framed surface in Definition \ref{involute-framed-surface} and the constant $c=0$.
\par
$(2)$ By a direct calculation, we have ${\bt}^{s,n}(u,v)=\sin \theta^{s,n}(u,v) \bt(u,v)+\cos \theta^{s,n}(u,v) \bn(u,v)$ and ${\bt}^{t,n}(u,v)=\sin\theta^{t,n}(u,v) \bn(u,v)+\cos\theta^{t,n}(u,v) \bs(u,v)$.
}
\end{remark}
\begin{corollary}\label{caustic_involute_basic_invariant}
Under the same notations as in Definitions \ref{caustic-framed-surface} and \ref{involute-framed-surface}, we have the following.
\par
$(1)$ The basic invariants of $\mathcal{C}^{s}(\bx,\bn,\bs)=({\bx}^{n,s},{\bn}^{n,s},{\bs}^{n,s})$ are given by 
\begin{align*}
&\begin{pmatrix}
{a}_1^{n,s} & {b}_1^{n,s}\\
{a}_2^{n,s} & {b}_2^{n,s}
\end{pmatrix}
=
\begin{pmatrix}
\lambda^{n,s}_u & (a_1+\lambda^{n,s} e_1) \cos \theta^{n,s}-(b_1+\lambda^{n,s} f_1) \sin \theta^{n,s}\\
\lambda^{n,s}_v & (a_2+\lambda^{n,s} e_2) \cos \theta^{n,s}-(b_2+\lambda^{n,s} f_2) \sin \theta^{n,s}
\end{pmatrix},\\
&\begin{pmatrix}
{e}_1^{n,s} & {f}_1^{n,s} & {g}_1^{n,s}\\
{e}_2^{n,s} & {f}_2^{n,s} & {g}_2^{n,s}
\end{pmatrix}
\\=&
\begin{pmatrix}
-e_1 \sin \theta^{n,s}-f_1\cos \theta^{n,s} & \theta^{n,s}_u-g_1 & e_1\cos \theta^{n,s}-f_1\sin \theta^{n,s}\\
-e_2 \sin \theta^{n,s}-f_2\cos \theta^{n,s} & \theta^{n,s}_v-g_2 & e_2\cos \theta^{n,s}-f_2\sin \theta^{n,s}
\end{pmatrix}.
\end{align*}
$(2)$ The basic invariants of $\mathcal{C}^{t}(\bx,\bn,\bs)=({\bx}^{n,t},{\bn}^{n,t},{\bs}^{n,t})$ are given by 
\begin{align*}
&\begin{pmatrix}
{a}_1^{n,t} & {b}_1^{n,t}\\
{a}_2^{n,t} & {b}_2^{n,t}
\end{pmatrix}
=
\begin{pmatrix}
(a_1+\lambda^{n,t} e_1) \sin\theta^{n,t}+(b_1+\lambda^{n,t} f_1) \cos\theta^{n,t} & \lambda^{n,t}_u\\
(a_2+\lambda^{n,t} e_2) \sin\theta^{n,t}+(b_2+\lambda^{n,t} f_2) \cos\theta^{n,t} & \lambda^{n,t}_v
\end{pmatrix},\\
&\begin{pmatrix}
{e}_1^{n,t} & {f}_1^{n,t} & {g}_1^{n,t}\\
{e}_2^{n,t} & {f}_2^{n,t} & {g}_2^{n,t}
\end{pmatrix}
\\=&
\begin{pmatrix}
g_1-\theta^{n,t}_u & -e_1 \cos\theta^{n,t}+f_1 \sin\theta^{n,t} & -e_1 \sin\theta^{n,t}-f_1 \cos\theta^{n,t}\\
g_2-\theta^{n,t}_v & -e_2 \cos\theta^{n,t}+f_2 \sin\theta^{n,t} & -e_2 \sin\theta^{n,t}-f_2 \cos\theta^{n,t}
\end{pmatrix}.
\end{align*}
$(3)$ The basic invariants of $\mathcal{I}^{s}(\bx,\bn,\bs)=({\bx}^{s,n},{\bn}^{s,n},{\bs}^{s,n})$ are given by 
\begin{align*}
&\begin{pmatrix}
{a}_1^{s,n} & {b}_1^{s,n}\\
{a}_2^{s,n} & {b}_2^{s,n}
\end{pmatrix}
\\=&
\begin{pmatrix}
(b_1+\lambda^{s,n} g_1) \cos\theta^{s,n}+\lambda^{s,n} e_1 \sin \theta^{s,n} &
(b_1+\lambda^{s,n} g_1) \sin\theta^{s,n}-\lambda^{s,n} e_1 \cos \theta^{s,n}\\
(b_2+\lambda^{s,n} g_2) \cos\theta^{s,n}+\lambda^{s,n} e_2 \sin \theta^{s,n} & 
(b_2+\lambda^{s,n} g_2) \sin\theta^{s,n}-\lambda^{s,n} e_2 \cos \theta^{s,n}
\end{pmatrix},\\
&\begin{pmatrix}
{e}_1^{s,n} & {f}_1^{s,n} & {g}_1^{s,n}\\
{e}_2^{s,n} & {f}_2^{s,n} & {g}_2^{s,n}
\end{pmatrix}
\\=&
\begin{pmatrix}
e_1 \sin\theta^{s,n}+g_1 \cos\theta^{s,n} & -e_1 \cos\theta^{s,n}+g_1 \sin\theta^{s,n} & -f_1-\theta^{s,n}_u\\
e_2 \sin\theta^{s,n}+g_2 \cos\theta^{s,n} & -e_2 \cos\theta^{s,n}+g_2 \sin\theta^{s,n} & -f_2-\theta^{s,n}_v
\end{pmatrix}.
\end{align*}
$(4)$ The basic invariants of $\mathcal{I}^{t}(\bx,\bn,\bs)=({\bx}^{t,n},{\bn}^{t,n},{\bs}^{t,n})$ are given by 
{\small \begin{align*}
&\begin{pmatrix}
{a}_1^{t,n} & {b}_1^{t,n}\\
{a}_2^{t,n} & {b}_2^{t,n}
\end{pmatrix}
\\=&
\begin{pmatrix}
-(a_1-\lambda^{t,n} g_1) \sin\theta^{t,n}-\lambda^{t,n} f_1 \cos\theta^{t,n} & (a_1-\lambda^{t,n} g_1) \cos\theta^{t,n}-\lambda^{t,n} f_1 \sin\theta^{t,n}\\
-(a_2-\lambda^{t,n} g_2) \sin\theta^{t,n}-\lambda^{t,n} f_2 \cos\theta^{t,n} & (a_2-\lambda^{t,n} g_2) \cos\theta^{t,n}-\lambda^{t,n} f_2 \sin\theta^{t,n}
\end{pmatrix},\\
&\begin{pmatrix}
{e}_1^{t,n} & {f}_1^{t,n} & {g}_1^{t,n}\\
{e}_2^{t,n} & {f}_2^{t,n} & {g}_2^{t,n}
\end{pmatrix}
\\=&
\begin{pmatrix}
-f_1 \cos\theta^{t,n}+g_1 \sin\theta^{t,n} & -f_1 \sin\theta^{t,n}-g_1 \cos\theta^{t,n} & e_1-\theta^{t,n}_u\\
-f_2 \cos\theta^{t,n}+g_2 \sin\theta^{t,n} & -f_2 \sin\theta^{t,n}-g_2 \cos\theta^{t,n} & e_2-\theta^{t,n}_v
\end{pmatrix}.
\end{align*}}
\end{corollary}
We consider conditions that caustics and involutes are inverse operations of framed surfaces.
\begin{theorem}\label{caustic-involute}
Let $(\bx,\bn,\bs):U \to \R^3 \times \Delta$ be a framed surface with basic invariants $(\mathcal{G},\mathcal{F}_1,\mathcal{F}_2)$. 
\par
$(1)$ $\rm(i)$ Suppose that ${\rm det}(\bb(u,v),\bg(u,v))=0$ for all $(u,v) \in U$, $\theta^{s,n}:U \to \R$ is a smooth function and a smooth function $\lambda^{s,n}:U \to \R$ is given by
$$
\lambda^{s,n}(u,v)=-\left(\int^u_{u_0} a_1(u,v) du+\int^v_{v_0} a_2(u_0,v) dv\right),
$$
for a point $(u_0,v_0) \in U$. If we take $\lambda^{n,s}, \theta^{n,s} : U \to \R$ by $\lambda^{n,s}=-\lambda^{s,n}$ and $\theta^{n,s}=-\theta^{s,n}$, then ${\mathcal{C}}^{s}({\mathcal{I}}^{s}(\bx,\bn,\bs))=({\bx},\bn,\bs)$.
\par
$\rm(ii)$ Suppose that ${\rm det}(\ba(u,v),\bg(u,v))=0$ for all $(u,v) \in U$, $\theta^{t,n}:U \to \R$ is a smooth function and a smooth function  $\lambda^{t,n}:U \to \R$ is given by
$$
\lambda^{t,n}(u,v)=-\left(\int^u_{u_0} b_1(u,v) du+\int^v_{v_0} b_2(u_0,v) dv\right),
$$
for a point $(u_0,v_0) \in U$. If we take $\lambda^{n,t}, \theta^{n,t} : U \to \R$ by $\lambda^{n,t}=-\lambda^{t,n}$ and $\theta^{n,t}=-\theta^{t,n}$, then ${\mathcal{C}}^{t}({\mathcal{I}}^{t}(\bx,\bn,\bs))=({\bx},\bn,\bs)$.
\par
$(2)$ $\rm(i)$ Suppose that there exist smooth functions $\lambda^{n,s}, \theta^{n,s}:U \to \R$ such that the condition \eqref{ns-condition2} satisfies. If we take $\theta^{s,n} : U \to \R$ by $\theta^{s,n}=-\theta^{n,s}$, then ${\mathcal{I}}^{s}({\mathcal{C}}^{s}(\bx,\bn,\bs))=(\bx+\lambda^{n,s}(u_0,v_0){\bn},\bn,\bs)$ for a point $(u_0,v_0) \in U$.
\par
$\rm(ii)$ Suppose that there exist smooth functions $\lambda^{n,t}, \theta^{n,t}:U \to \R$ such that the condition \eqref{nt-condition2} satisfies. If we take $\theta^{t,n} : U \to \R$ by $\theta^{t,n}=-\theta^{n,t}$, then ${\mathcal{I}}^{t}({\mathcal{C}}^{t}(\bx,\bn,\bs))=(\bx+\lambda^{n,t}(u_0,v_0){\bn},\bn,\bs)$ for a point $(u_0,v_0) \in U$.
\end{theorem}
\demo
$(1)$ $\rm(i)$ By Definition \ref{involute-framed-surface} $(1)$, the map $\mathcal{I}^{s}:\mathcal{F}(U,\R^3 \times \Delta)\to\mathcal{F}(U,\R^3 \times \Delta)$ is given by $\mathcal{I}^{s}(\bx,\bn,\bs)=({\bx}^{s,n},{\bn}^{s,n},{\bs}^{s,n})=({\bx}+\lambda^{s,n}\bs,\bs,\cos\theta^{s,n}\bt+\sin\theta^{s,n} \bn)$.
By Corollary \ref{caustic_involute_basic_invariant} (3), the basic invariants of $\mathcal{I}^{s}(\bx,\bn,\bs)$ is given by $({\mathcal{G}}^{s,n},{\mathcal{F}}_1^{s,n},{\mathcal{F}}_2^{s,n})$. 
The condition (\ref{ns-condition2}) for $\mathcal{I}^{s}(\bx,\bn,\bs)$ is given by
\begin{align*}
\begin{pmatrix}
a_1^{s,n}+\lambda^{n,s} e_1^{s,n} & b_1^{s,n}+\lambda^{n,s} f_1^{s,n}\\
a_2^{s,n}+\lambda^{n,s} e_2^{s,n} & b_2^{s,n}+\lambda^{n,s} f_2^{s,n}
\end{pmatrix}
\begin{pmatrix}
\sin\theta^{n,s} \\
\cos\theta^{n,s}
\end{pmatrix}
=
\begin{pmatrix}
0 \\
0
\end{pmatrix}.
\end{align*}
By a direct calculation, we have
\begin{align}\label{CsIs-condition}
\begin{pmatrix}
b_1+(\lambda^{n,s}+\lambda^{s,n}) g_1 & (\lambda^{n,s}+\lambda^{s,n}) e_1\\
b_2+(\lambda^{n,s}+\lambda^{s,n}) g_2 & (\lambda^{n,s}+\lambda^{s,n}) e_2
\end{pmatrix}
\begin{pmatrix}
\sin(\theta^{n,s}+\theta^{s,n}) \\
-\cos(\theta^{n,s}+\theta^{s,n})
\end{pmatrix}
=
\begin{pmatrix}
0 \\
0
\end{pmatrix}.
\end{align}
If we take $\lambda^{n,s}, \theta^{n,s} : U \to \R$ by $\lambda^{n,s}=-\lambda^{s,n}$ and $\theta^{n,s}=-\theta^{s,n}$, then the condition (\ref{CsIs-condition}) is satisfied. 
Thus, the map ${\mathcal{C}}^{s}$ of the map $\mathcal{I}^{s}$ exists.
By Definition \ref{caustic-framed-surface} $(1)$, the map ${\mathcal{C}}^{s}$ of the map $\mathcal{I}^{s}$, 
$$
\mathcal{C}^{s}(\mathcal{I}^{s}(\bx,\bn,\bs))
=({\bx}^{n,s}({\bx}^{s,n}, {\bn}^{s,n}, {\bs}^{s,n}),{\bn}^{n,s}({\bx}^{s,n}, {\bn}^{s,n}, {\bs}^{s,n}),{\bs}^{n,s}({\bx}^{s,n}, {\bn}^{s,n}, {\bs}^{s,n}))
$$
is given by 
\begin{align*}
{\bx}^{n,s}({\bx}^{s,n}, {\bn}^{s,n}, {\bs}^{s,n})&={\bx}^{s,n}+{\lambda^{n,s}}{\bn}^{s,n}=\bx+(\lambda^{n,s}+\lambda^{s,n})\bs=\bx,\\
{\bn}^{n,s}({\bx}^{s,n}, {\bn}^{s,n}, {\bs}^{s,n})&=\sin\theta^{n,s}{\bs}^{s,n}+\cos\theta^{n,s}{\bt}^{s,n}\\
&=(-\sin\theta^{n,s}\sin\theta^{s,n}+\cos\theta^{n,s}\cos\theta^{s,n})\bn\\
&\qquad+(\sin\theta^{n,s}\cos\theta^{s,n}
+\cos\theta^{n,s}\sin\theta^{s,n})\bt\\
&=\cos(\theta^{n,s}+\theta^{s,n})\bn+\sin(\theta^{n,s}+\theta^{s,n})\bt=\bn,\\
{\bs}^{n,s}({\bx}^{s,n}, {\bn}^{s,n}, {\bs}^{s,n})&={\bn}^{s,n}=\bs.
\end{align*}\par
$\rm(ii)$ We can also prove by the same method of $\rm(i)$.
\par
$(2)$ $\rm(i)$ By Definition \ref{caustic-framed-surface} $(1)$, the map $\mathcal{C}^{s}:\mathcal{F}(U,\R^3 \times \Delta)\to\mathcal{F}(U,\R^3 \times \Delta)$ is given by $\mathcal{C}^{s}(\bx,\bn,\bs)=({\bx}^{n,s},{\bn}^{n,s},{\bs}^{n,s})=({\bx}+\lambda^{n,s}\bn,\sin \theta^{n,s}\bs+\cos\theta^{n,s}\bt,\bn)$,  where there exist smooth functions $\lambda^{n,s}, \theta^{n,s}:U \to \R$ such that the condition (\ref{ns-condition2}) satisfies. By corollary \ref{caustic_involute_basic_invariant} (1), the basic invariants of $\mathcal{C}^s(\bx,\bn,\bs)$ is given by $({\mathcal{G}}^{n,s},{\mathcal{F}}_1^{n,s},{\mathcal{F}}_2^{n,s})$. 
By the integrability conditions (\ref{integrability.condition}), we have
\begin{align*}
{\rm det}(\bb^{n,s}(u,v),\bg^{n,s}(u,v))&={{a}_{1v}^{n,s}}(u,v)-{{a}_{2u}^{n,s}}(u,v)={{\lambda}_{uv}^{n,s}}(u,v)-{{\lambda}_{vu}^{n,s}}(u,v)=0,
\end{align*}
for all $(u,v) \in U$. Thus, the map $\mathcal{I}^{s}$ of the map $\mathcal{C}^{s}$ always exists.
By Definition \ref{involute-framed-surface} $(1)$, $\lambda^{s,n}:U\to\R$ is given by
\begin{align*}
\lambda^{s,n}(u,v)&=-\left(\int^u_{u_0} {a}_1^{n,s}(u,v) du+\int^v_{v_0} {a}_2^{n,s}(u_0,v) dv \right)\\
&=-\left(\int^u_{u_0} \lambda^{n,s}_u(u,v) du+\int^v_{v_0} \lambda^{n,s}_v(u_0,v) dv \right)\\
&=-\left(\lambda^{n,s}(u,v)-\lambda^{n,s}(u_0,v_0)\right),
\end{align*}
for a point $(u_0,v_0) \in U$. if we take $\theta^{s,n} : U \to \R$ by $\theta^{s,n}=-\theta^{n,s}$, the map ${\mathcal{I}}^{s}$ of the map $\mathcal{C}^{s}$, 
$$
\mathcal{I}^{s}(\mathcal{C}^{s}(\bx,\bn,\bs))
=({\bx}^{s,n}({\bx}^{n,s}, {\bn}^{n,s}, {\bs}^{n,s}), {\bn}^{s,n}({\bx}^{n,s}, {\bn}^{n,s}, {\bs}^{n,s}), {\bs}^{s,n}({\bx}^{n,s}, {\bn}^{n,s}, {\bs}^{n,s}))
$$
is given by 
\begin{align*}
{\bx}^{s,n}({\bx}^{n,s}, {\bn}^{n,s}, {\bs}^{n,s})&={\bx}^{n,s}+{\lambda^{s,n}}{\bs}^{n,s}=\bx+(\lambda^{n,s}+\lambda^{s,n})\bn\\
&=\bx+\lambda^{n,s}(u_0,v_0)\bn,\\
{\bn}^{s,n}({\bx}^{n,s}, {\bn}^{n,s}, {\bs}^{n,s})&={\bs}^{n,s}=\bn,\\
{\bs}^{s,n}({\bx}^{n,s}, {\bn}^{n,s}, {\bs}^{n,s})&=\cos\theta^{s,n}{\bt}^{n,s}-\sin\theta^{s,n}{\bn}^{n,s}\\
&=(-\sin\theta^{n,s}\sin\theta^{s,n}
+\cos\theta^{n,s}\cos\theta^{s,n})\bs\\
&\quad-(\sin\theta^{n,s}\cos\theta^{s,n}
+\cos\theta^{n,s}\sin\theta^{s,n})\bt\\
&=\cos(\theta^{n,s}+\theta^{s,n})\bs-\sin(\theta^{n,s}+\theta^{s,n})\bt=\bs.
\end{align*}\par
$\rm(ii)$ We can also prove by the same method of $\rm(i)$.
\enD
\begin{theorem}\label{caustic_involute_to_SS}
Let $(\bx,\bn,\bs):U \to \R^3 \times \Delta$ be a framed surface with basic invariants $(\mathcal{G},\mathcal{F}_1,\mathcal{F}_2)$. 
\par
$(1)$ Suppose that ${\rm det}(\bb(u,v),\bg(u,v))=0$ for all $(u,v) \in U$, $\theta^{s,n}:U \to \R$ is a smooth function and a smooth function  $\lambda^{s,n}:U \to \R$ is given by
$$
\lambda^{s,n}(u,v)=-\left(\int^u_{u_0} a_1(u,v) du+\int^v_{v_0} a_2(u_0,v) dv\right),
$$
for a point $(u_0,v_0) \in U$. 
Then $(\bx,\bn,\bs)$ is an $(\bs,\overline{\bs})$-Bertrand framed surface if and only if there exists a function $\lambda^{n,s}:U \to \R$ with $\lambda^{n,s}+\lambda^{s,n} \not\equiv 0$ such that the map ${\mathcal{C}}^{s}$ of the map $\mathcal{I}^{s}$ exists. 
Moreover, $(\bx,\bn,\bs)$ and ${\mathcal{C}}^{s}({\mathcal{I}}^{s}(\bx,\bn,\bs))$  are $(\bs, \overline{\bs})$-mates.
\par
$(2)$ Suppose that ${\rm det}(\ba(u,v),\bg(u,v))=0$ for all $(u,v) \in U$, $\theta^{t,n}:U \to \R$ is a smooth function and a smooth function  $\lambda^{t,n}:U \to \R$ is given by
$$
\lambda^{t,n}(u,v)=-\left(\int^u_{u_0} b_1(u,v) du+\int^v_{v_0} b_2(u_0,v) dv\right),
$$
for a point $(u_0,v_0) \in U$. 
Then $(\bx,\bn,\bs)$ is a $(\bt,\overline{\bt})$-Bertrand framed surface if and only if there exists a function $\lambda^{n,t}:U \to \R$ with $\lambda^{n,t}+\lambda^{t,n} \not\equiv 0$ such that the map the map ${\mathcal{C}}^{t}$ of the map $\mathcal{I}^{t}$ exists. 
Moreover, $(\bx,\bn,\bs)$ and ${\mathcal{C}}^{t}({\mathcal{I}}^{t}(\bx,\bn,\bs))$ are $(\bt, \overline{\bt})$-mates.
\end{theorem}
\demo
$(1)$ Suppose that $(\bx,\bn,\bs)$ is an $(\bs,\overline{\bs})$-Bertrand framed surface. 
By Definition \ref{involute-framed-surface} $(1)$, the map $\mathcal{I}^{s}:\mathcal{F}(U,\R^3 \times \Delta)\to\mathcal{F}(U,\R^3 \times \Delta)$ is given by $\mathcal{I}^{s}(\bx,\bn,\bs)=({\bx}^{s,n},{\bn}^{s,n},{\bs}^{s,n})=({\bx}+\lambda^{s,n}\bs,\bs,\cos\theta^{s,n}\bt+\sin\theta^{s,n} \bn)$.
By Corollary \ref{caustic_involute_basic_invariant} (3), the basic invariants of $\mathcal{I}^{s}(\bx,\bn,\bs)$ is given by $({\mathcal{G}}^{s,n},{\mathcal{F}}_1^{s,n},{\mathcal{F}}_2^{s,n})$. 
The condition (\ref{ns-condition2}) for $\mathcal{I}^{s}(\bx,\bn,\bs)$ is given by the condition (\ref{CsIs-condition}). 
By 
Theorem \ref{ss-Bertrand}, $(\overline{\bx},\overline{\bn},\overline{\bs}):U \to \R^3 \times \Delta$ is given by $(\overline{\bx},\overline{\bn},\overline{\bs})=(\bx+\lambda \bs, \sin \theta \bt+\cos \theta \bn,\bs)$, where there exist smooth functions $\lambda, \theta:U \to \R$ with $\lambda \not\equiv 0$ such that the condition \eqref{ss-condition} satisfies. 
If we take $\lambda^{n,s}, \theta^{n,s}:U\to\R$ by $\lambda^{n,s}=\lambda-\lambda^{s,n}$ and $\theta^{n,s}=\theta-\theta^{s,n}$, then the condition (\ref{CsIs-condition}) is satisfied. Thus, the map ${\mathcal{C}}^{s}$ of the map $\mathcal{I}^{s}$ exists.
By Definition \ref{caustic-framed-surface} $(1)$, the map ${\mathcal{C}}^{s}$ of the map $\mathcal{I}^{s}$,
$$
\mathcal{C}^{s}(\mathcal{I}^{s}(\bx,\bn,\bs))
=({\bx}^{n,s}({\bx}^{s,n}, {\bn}^{s,n}, {\bs}^{s,n}),{\bn}^{n,s}({\bx}^{s,n}, {\bn}^{s,n}, {\bs}^{s,n}),{\bs}^{n,s}({\bx}^{s,n}, {\bn}^{s,n}, {\bs}^{s,n}))
$$
is given by 
\begin{align*}
{\bx}^{n,s}({\bx}^{s,n}, {\bn}^{s,n}, {\bs}^{s,n})&={\bx}^{s,n}+{\lambda^{n,s}}{\bn}^{s,n}=\bx+(\lambda^{n,s}+\lambda^{s,n})\bs=\bx+\lambda\bs=\overline{\bx},\\
{\bn}^{n,s}({\bx}^{s,n}, {\bn}^{s,n}, {\bs}^{s,n})&=\sin\theta^{n,s}{\bs}^{s,n}+\cos\theta^{n,s}{\bt}^{s,n}\\
&=(-\sin\theta^{n,s}\sin\theta^{s,n}+\cos\theta^{n,s}\cos\theta^{s,n})\bn\\
&\qquad +(\sin\theta^{n,s}\cos\theta^{s,n}+\cos\theta^{n,s}\sin\theta^{s,n})\bt\\
&=\cos(\theta^{n,s}+\theta^{s,n})\bn+\sin(\theta^{n,s}+\theta^{s,n})\bt=\cos\theta\bn+\sin\theta\bt=\overline{\bn},\\
{\bs}^{n,s}({\bx}^{s,n}, {\bn}^{s,n}, {\bs}^{s,n})&={\bn}^{s,n}=\bs=\overline{\bs}. 
\end{align*}
Therefore, $(\bx,\bn,\bs)$ and ${\mathcal{C}}^{s}({\mathcal{I}}^{s}(\bx,\bn,\bs))$ are $(\bs, \overline{\bs})$-mates.
\par 
Conversely, by Definition \ref{involute-framed-surface} $(1)$, the map $\mathcal{I}^{s}:\mathcal{F}(U,\R^3 \times \Delta)\to\mathcal{F}(U,\R^3 \times \Delta)$ is given by $\mathcal{I}^{s}(\bx,\bn,\bs)=({\bx}^{s,n},{\bn}^{s,n},{\bs}^{s,n})=({\bx}+\lambda^{s,n}\bs,\bs,\cos\theta^{s,n}\bt+\sin\theta^{s,n} \bn)$.
By Corollary \ref{caustic_involute_basic_invariant} (3), the basic invariants of $\mathcal{I}^{s}(\bx,\bn,\bs)$ is given by $({\mathcal{G}}^{s,n},{\mathcal{F}}_1^{s,n},{\mathcal{F}}_2^{s,n})$. 
By assumption, there exist smooth functions $\lambda^{n,s}, \theta^{n,s}:U \to \R$ with $\lambda^{n,s}+\lambda^{s,n}\neq 0$ and the condition (\ref{CsIs-condition}) satisfies. If we take $\lambda, \theta:U \to \R$ by $\lambda=\lambda^{n,s}+\lambda^{s,n}$ and $\theta=\theta^{n,s}+\theta^{s,n}$, then $\lambda \not\equiv 0$ and the condition (\ref{ss-condition}) is satisfied. By Theorem \ref{ss-Bertrand}, $(\bx,\bn,\bs)$ is an $(\bs, \overline{\bs})$-Bertrand framed surface. By the same calculations in the proof of the Theorem \ref{caustic_involute_to_SS} (1), $(\bx,\bn,\bs)$ and ${\mathcal{C}}^{s}({\mathcal{I}}^{s}(\bx,\bn,\bs))$ are $(\bs, \overline{\bs})$-mates.
\par
$(2)$ We can also prove by the same method of $(1)$.
\enD
\section{Tangential direction framed surfaces}

Let $(\bx,\bn,\bs):U \to \R^3 \times \Delta$ be a framed surface with basic invariants $(\mathcal{G},\mathcal{F}_1,\mathcal{F}_2)$. 
\begin{definition}\label{tangent-framed-surface}{\rm 
$(1)$
The map $\mathcal{S}^{t}:\mathcal{F}(U,\R^3 \times \Delta)\to\mathcal{F}(U,\R^3 \times \Delta)$, $\mathcal{S}^{t}(\bx,\bn,\bs)=({\bx}^{s,t},{\bn}^{s,t},{\bs}^{s,t})$ is given by
\begin{align*}
{\bx}^{s,t}(u,v)&={\bx}(u,v)+{\lambda^{s,t}(u,v)}{\bs}(u,v),\\
{\bn}^{s,t}(u,v)&=\cos\theta^{s,t}(u,v)\bt(u,v)-\sin\theta^{s,t}(u,v)\bn(u,v),\\
{\bs}^{s,t}(u,v)&=\sin\theta^{s,t}(u,v)\bt(u,v)+\cos\theta^{s,t}(u,v)\bn(u,v),
\end{align*}
where there exist smooth functions $\lambda^{s,t}, \theta^{s,t}:U \to \R$ such that 
\begin{eqnarray}\label{st-condition2}
\begin{pmatrix}
\lambda^{s,t}(u,v)e_1(u,v) & b_1(u,v)+\lambda^{s,t}(u,v)g_1(u,v)\\
\lambda^{s,t}(u,v)e_2(u,v) & b_2(u,v)+\lambda^{s,t}(u,v)g_2(u,v)
\end{pmatrix}
\begin{pmatrix}
\sin {\theta^{s,t}}(u,v)\\
\cos {\theta^{s,t}}(u,v)
\end{pmatrix}
=
\begin{pmatrix}
0\\
0
\end{pmatrix}
\end{eqnarray}
for all $(u,v) \in U$. 
We say that $({\bx}^{s,t},{\bn}^{s,t},{\bs}^{s,t})$ is a {\it tangential direction framed surface with respect to $\bs$} of the framed surface $(\bx,\bn,\bs)$. \par
$(2)$
The map $\mathcal{T}^{s}:\mathcal{F}(U,\R^3 \times \Delta)\to\mathcal{F}(U,\R^3 \times \Delta)$, $\mathcal{T}^{s}(\bx,\bn,\bs)=({\bx}^{t,s},{\bn}^{t,s},{\bs}^{t,s})$ is given by
\begin{align*}
{\bx}^{t,s}(u,v)&={\bx}(u,v)+{\lambda^{t,s}(u,v)}{\bt}(u,v),\\
{\bn}^{t,s}(u,v)&=\sin\theta^{t,s}(u,v)\bn(u,v)+\cos\theta^{t,s}(u,v)\bs(u,v),\\
{\bs}^{t,s}(u,v)&=\bt(u,v),
\end{align*}
where there exist smooth functions $\lambda^{t,s}, \theta^{t,s}:U \to \R$ such that 
\begin{eqnarray}\label{ts-condition2}
\begin{pmatrix}
\lambda^{t,s}(u,v)f_1(u,v) & a_1(u,v)-\lambda^{t,s}(u,v)g_1(u,v)\\
\lambda^{t,s}(u,v)f_2(u,v) & a_2(u,v)-\lambda^{t,s}(u,v)g_2(u,v)
\end{pmatrix}
\begin{pmatrix}
-\sin {\theta^{t,s}}(u,v)\\
\cos {\theta^{t,s}}(u,v)
\end{pmatrix}
=
\begin{pmatrix}
0\\
0
\end{pmatrix}
\end{eqnarray}
for all $(u,v) \in U$. 
Then we say that $({\bx}^{t,s},{\bn}^{t,s},{\bs}^{t,s})$ is a {\it tangential direction framed surface with respect to $\bt$} of the framed surface $(\bx,\bn,\bs)$. 
}
\end{definition}
\begin{remark}{\rm
$(1)$ The map $\mathcal{S}^{t}$ (respectively, $\mathcal{T}^{s}$) is corresponding to the $(\bs, \overline{\bt})$ (respectively, $(\bt, \overline{\bs})$)-Bertrand framed surface. 
\par
$(2)$ By a direct calculation, we have ${\bt}^{s,t}(u,v)=\bs(u,v)$ and ${\bt}^{s,t}(u,v)=\cos\theta^{t,s}(u,v)\bn(u,v)-\sin\theta^{t,s}(u,v)\bs(u,v)$.
}
\end{remark}
\begin{corollary}\label{StTs_basic_invariant}{ 
Under the same notations as in Definition \ref{tangent-framed-surface}, we have the following.
\begin{itemize}
\item[$(1)$] The basic invariants of $\mathcal{S}^{t}(\bx,\bn,\bs)=({\bx}^{s,t},{\bn}^{s,t},{\bs}^{s,t})$ are given by 
\begin{align*}
&\begin{pmatrix}
{a}_1^{s,t} & {b}_1^{s,t}\\
{a}_2^{s,t} & {b}_2^{s,t}
\end{pmatrix}
=
\begin{pmatrix}
-\lambda^{s,t} e_1 \cos \theta^{s,t}+(b_1+\lambda^{s,t} g_1) \sin \theta^{s,t} & a_1+\lambda^{s,t}_u\\
-\lambda^{s,t} e_2 \cos \theta^{s,t}+(b_2+\lambda^{s,t} g_2) \sin \theta^{s,t} & a_2+\lambda^{s,t}_v 
\end{pmatrix},\\
&\begin{pmatrix}
{e}_1^{s,t} & {f}_1^{s,t} & {g}_1^{s,t}\\
{e}_2^{s,t} & {f}_2^{s,t} & {g}_2^{s,t}
\end{pmatrix}\\
=&
\begin{pmatrix}
-\theta^{s,t}_u-f_1 & -g_1\cos\theta^{s,t}-e_1\sin\theta^{s,t} & -g_1\sin\theta^{s,t}+e_1\cos\theta^{s,t} \\
-\theta^{s,t}_v-f_2 & -g_2\cos\theta^{s,t}-e_2\sin\theta^{s,t} & -g_2\sin\theta^{s,t}+e_2\cos\theta^{s,t}
\end{pmatrix}.
\end{align*} 
\item[$(2)$] The basic invariants of $\mathcal{T}^{s}(\bx,\bn,\bs)=({\bx}^{t,s},{\bn}^{t,s},{\bs}^{t,s})$ are given by 
\begin{align*}
&\begin{pmatrix}
{a}_1^{t,s} & {b}_1^{t,s}\\
{a}_2^{t,s} & {b}_2^{t,s}
\end{pmatrix}=
\begin{pmatrix}
b_1+\lambda^{t,s}_u & -\lambda^{t,s} f_1 \cos \theta^{t,s}-(a_1-\lambda^{t,s} g_1) \sin \theta^{t,s}\\
b_2+\lambda^{t,s}_v & -\lambda^{t,s} f_2 \cos \theta^{t,s}-(a_2-\lambda^{t,s} g_2) \sin \theta^{t,s}
\end{pmatrix},\\
&\begin{pmatrix}
{e}_1^{t,s} & {f}_1^{t,s} & {g}_1^{t,s}\\
{e}_2^{t,s} & {f}_2^{t,s} & {g}_2^{t,s}
\end{pmatrix}\\
=&
\begin{pmatrix}
f_1 \sin \theta^{t,s}+g_1\cos \theta^{t,s} & \theta^{t,s}_u-e_1 & -f_1\cos \theta^{t,s}+g_1\sin \theta^{t,s}\\
f_2 \sin \theta^{t,s}+g_2\cos \theta^{t,s} & \theta^{t,s}_v-e_2 & -f_2\cos \theta^{t,s}+g_2\sin \theta^{t,s}
\end{pmatrix}.
\end{align*}
\end{itemize}
}
\end{corollary}
We give conditions that tangential direction framed surfaces are inverse operations of framed surfaces.
\begin{theorem}\label{S-T}
Let $(\bx,\bn,\bs):U \to \R^3 \times \Delta$ be a framed surface. 
\par
$(1)$ Suppose that there exist smooth functions $\lambda^{s,t}, \theta^{s,t}:U \to \R$ such that the condition \eqref{st-condition2} satisfies. If we take $\lambda^{t,s}, \theta^{t,s} : U \to \R$ by $\lambda^{t,s}=-\lambda^{s,t}$ and $\theta^{t,s}=-\theta^{s,t}$, then ${\mathcal{T}}^{s}({\mathcal{S}}^{t}(\bx,\bn,\bs))=(\bx,\bn,\bs)$.
\par
$(2)$ Suppose that there exist smooth functions $\lambda^{t,s}, \theta^{t,s}:U \to \R$ such that the condition \eqref{ts-condition2} satisfies. If we take $\lambda^{s,t}, \theta^{s,t} : U \to \R$ by $\lambda^{s,t}=-\lambda^{t,s}$ and $\theta^{s,t}=-\theta^{t,s}$, then ${\mathcal{S}}^{t}({\mathcal{T}}^{s}(\bx,\bn,\bs))=(\bx,\bn,\bs)$.
\end{theorem}
\demo
$(1)$ By Definition \ref{tangent-framed-surface} $(1)$, the map $\mathcal{S}^{t}:\mathcal{F}(U,\R^3 \times \Delta)\to\mathcal{F}(U,\R^3 \times \Delta)$ is given by $\mathcal{S}^{t}(\bx,\bn,\bs)=({\bx}^{s,t},{\bn}^{s,t},{\bs}^{s,t})=({\bx}+\lambda^{s,t}\bs,\cos \theta^{s,t}\bt-\sin\theta^{s,t}\bn,\sin\theta^{s,t}\bt+\cos\theta^{s,t}\bn)$ where 
there exist smooth functions $\lambda^{s,t}, \theta^{s,t}:U \to \R$ such that the condition (\ref{st-condition2}) satisfies. By Corollary \ref{StTs_basic_invariant} (1), the basic invariants of $\mathcal{S}^{t}(\bx,\bn,\bs)$ is given by $({\mathcal{G}}^{\mathcal{S},t},{\mathcal{F}}_1^{\mathcal{S},t},{\mathcal{F}}_2^{\mathcal{S},t})$.
The condition (\ref{ts-condition2}) for $\mathcal{S}^{t}(\bx,\bn,\bs)$ is 
\begin{eqnarray}\label{StTs-condition}
\begin{pmatrix}
\lambda^{t,s}f_1^{s,t} & a_1^{s,t}-\lambda^{t,s}g_1^{s,t}\\
\lambda^{t,s}f_2^{s,t} & a_2^{s,t}-\lambda^{t,s}g_2^{s,t}
\end{pmatrix}
\begin{pmatrix}
-\sin {\theta^{t,s}}\\
\cos {\theta^{t,s}}
\end{pmatrix}
=
\begin{pmatrix}
0\\
0
\end{pmatrix}.
\end{eqnarray}
If we take $\lambda^{t,s}, \theta^{t,s} : U \to \R$ by $\lambda^{t,s}=-\lambda^{s,t}$ and $\theta^{t,s}=-\theta^{s,t}$, we have 
\begin{align*}
&-\lambda^{t,s}f_i^{s,t}\sin {\theta^{t,s}}+(a_i^{s,t}-\lambda^{t,s}g_i^{s,t})\cos{\theta^{t,s}}\\
&=-\lambda^{t,s}(-g_i\cos\theta^{s,t}-e_i\sin\theta^{s,t})\\
&\quad+\left(-\lambda^{s,t}e_i\cos\theta^{s,t}+(b_i+\lambda^{s,t}g_i)\sin\theta^{s,t}-\lambda^{t,s}(-g_i\sin\theta^{s,t}+e_i\cos\theta^{s,t})\right)\cos\theta^{t,s}\\
&=\lambda^{t,s}g_i\sin(\theta^{s,t}+\theta^{t,s})-\lambda^{t,s}e_i\cos(\theta^{s,t}+\theta^{t,s})\\
&\quad+\left(-\lambda^{s,t}e_i\cos\theta^{s,t}+(b_i+\lambda^{s,t}g_i)\sin\theta^{s,t}\right)\cos\theta^{t,s}\\
&=\lambda^{s,t}e_i-\lambda^{s,t}e_i\cos^2\theta^{s,t}+(b_i+\lambda^{s,t}g_i)\sin\theta^{s,t}\cos\theta^{s,t}\\
&=\lambda^{s,t}e_i\sin^2\theta^{s,t}+(b_i+\lambda^{s,t}g_i)\sin\theta^{s,t}\cos\theta^{s,t}\\
&=\sin\theta^{s,t}\left(\lambda^{s,t}e_i\sin\theta^{s,t}+(b_i+\lambda^{s,t}g_i)\cos\theta^{s,t}\right)\\
&=0,
\end{align*}
for $i=1,2$. It follows that the condition (\ref{StTs-condition}) is satisfied.
Thus, the map ${\mathcal{T}^{s}}$ of the map $\mathcal{S}^{t}$ exists.
By Definition \ref{tangent-framed-surface} $(2)$, the map ${\mathcal{T}^{s}}$ of the map $\mathcal{S}^{t}$,
$$
\mathcal{T}^{s}(\mathcal{S}^{t}(\bx,\bn,\bs))=({\bx}^{t,s}({\bx}^{s,t},{\bn}^{s,t},{\bs}^{s,t}),{\bn}^{t,s}({\bx}^{s,t},{\bn}^{s,t},{\bs}^{s,t}),{\bs}^{t,s}({\bx}^{s,t},{\bn}^{s,t},{\bs}^{s,t}))
$$
is given by 
\begin{align*}
{\bx}^{t,s}({\bx}^{s,t},{\bn}^{s,t},{\bs}^{s,t})&={\bx}^{s,t}+{\lambda^{t,s}}{\bs}^{t,s}=\bx+(\lambda^{s,t}+\lambda^{t,s})\bs=\bx,\\
{\bn}^{t,s}({\bx}^{s,t},{\bn}^{s,t},{\bs}^{s,t})&=\sin\theta^{t,s}{\bn}^{s,t}+\cos\theta^{t,s}{\bs}^{s,t}\\
&=(-\sin\theta^{s,t}\sin\theta^{t,s}+\cos\theta^{s,t}\cos\theta^{t,s})\bn\\
&\qquad-(\sin\theta^{s,t}\cos\theta^{t,s}+\cos\theta^{s,t}\sin\theta^{t,s})\bt\\
&=\cos(\theta^{s,t}+\theta^{t,s})\bn-\sin(\theta^{s,t}+\theta^{t,s})\bt=\bn,\\
{\bs}^{t,s}({\bx}^{s,t},{\bn}^{s,t},{\bs}^{s,t})&=\bt^{s,t}=\bs.
\end{align*}
\par
$(2)$ We can also prove by the same method of $(1)$.
\enD
\begin{remark}{\rm
If $e_1(u,v)=0$ and $e_2(u,v)=0$ (respectively, $f_1(u,v)=0$ and $f_2(u,v)=0$) for all $(u,v) \in U$, then $(\bx,\bn,\bs)$ is always an $(\bs,\overline{\bt})$ (respectively, $(\bt,\overline{\bs})$)-Bertrand framed surface for any $\lambda^{s,t} : U\to\R$ and for any constant $\theta^{s,t}$ with $\cos\theta^{s,t}=0$ (respectively, for any ${\lambda}^{t,s} : U\to\R$ and for any constant $\theta^{t,s}$ with $\cos\theta^{t,s}=0$).
}
\end{remark}
\section{Examples}

We give concrete examples of caustics, involutes and tangential direction framed surfaces.
\begin{example}{\rm(A cuspidal edge)
Let $(\bx,\bn,\bs) : \R^2\to \R^3 \times \Delta$ be 
$$
\bx(u,v)=\left(u,\frac{v^2}{2},\frac{v^3}{3}\right),\ \bn(u,v)=\frac{1}{\sqrt{v^2+1}}(0,-v,1), \ \bs(u,v)=(1,0,0).
$$
Then $\bt(u,v)=(0,1,v)/\sqrt{1+v^2}$ and $(\bx,\bn,\bs)$ is a framed surface with the basic invariants
\begin{align*}
\begin{pmatrix}
{a}_1 & {b}_1\\
{a}_2 & {b}_2
\end{pmatrix}
=
\begin{pmatrix}
1 & 0\\
0 & v\sqrt{v^2+1}
\end{pmatrix}, \ 
\begin{pmatrix}
{e}_1 & {f}_1 & {g}_1\\
{e}_2 & {f}_2 & {g}_2
\end{pmatrix}
=
\begin{pmatrix}
0 & 0 & 0\\
0 & -1/(v^2+1) & 0
\end{pmatrix}.
\end{align*}
It follows that the curvature $C^F$ of $(\bx, \bn, \bs)$ is given by
$$
J^F(u,v)=v\sqrt{v^2+1}, \ K^F(u,v)=0, \ H^F(u,v)=\frac{1}{2(v^2+1)}.
$$
If we take $\lambda^{n,s}(u,v)=v(v^2+1)^{3/2}$ and $\theta^{n,s}(u,v)=0$, then condition (\ref{ns-condition2}) is satisfied. 
Therefore, we have a caustic of the framed surface, $\mathcal{C}^{s}(\bx,\bn,\bs)=({\bx}^{n,s},{\bn}^{n,s},{\bs}^{n,s})$,
\begin{align*}
{\bx}^{n,s}(u,v)=\left(u, -v^4-\frac{v^2}{2}, \frac{4}{3}v^3+v\right), \ 
{\bn}^{n,s}(u,v)=\bt(u,v),\ {\bs}^{n,s}(u,v)= \bn(u,v).
\end{align*}
Moreover, if we take $\lambda^{n,t}(u,v)=v(v^2+1)^{3/2}$ and $\theta^{n,t}(u,v)=-\pi/2$, then condition (\ref{nt-condition2}) is satisfied. 
Therefore,we laso have a caustic of the framed surface, $\mathcal{C}^{t}(\bx,\bn,\bs)=({\bx}^{n,t},{\bn}^{n,t},{\bs}^{n,t})$, 
\begin{align*}
{\bx}^{n,t}(u,v)=\left(u, -v^4-\frac{v^2}{2}, \frac{4}{3}v^3+v\right), \ 
{\bn}^{n,t}(u,v)=\bt(u,v), \ {\bs}^{n,t}(u,v)=-\bs(u,v).
\end{align*}
Since ${\rm det}(\bb(u,v),\bg(u,v))=0$ for all $(u,v) \in U$, if we take $\lambda^{s,n}(u,v)=-u$, $\theta^{s,n}(u,v)=-\pi/2$ and $(u_0,v_0)=(0,0)$, then we have an involute with respect to $\bs$ at $(0,0)$, $\mathcal{I}^{s}(\bx,\bn,\bs)=({\bx}^{s,n},{\bn}^{s,n},{\bs}^{s,n}),$
\begin{align*}
{\bx}^{s,n}(u,v)=\left(0,\frac{v^2}{2},\frac{v^3}{3}\right), \ {\bn}^{s,n}(u,v)=\bs(u,v),\ {\bs}^{s,n}(u,v)=\bn(u,v).
\end{align*}
Moreover, since ${\rm det}(\ba(u,v),\bg(u,v))=0$ for all $(u,v) \in U$, if we take $\lambda^{t,n}(u,v)=-\frac{1}{3}\big((v^2+1)^{\frac{3}{2}}-1\big)$, $\theta^{t,n}(u,v)=0$ and $(u_0,v_0)=(0,0)$, then we have an involute with respect to $\bt$ at $(0,0)$, $\mathcal{I}^{t}(\bx,\bn,\bs)=({\bx}^{t,n},{\bn}^{t,n},{\bs}^{t,n}),$
\begin{align*}
&{\bx}^{t,n}(u,v)=\left(u,\frac{v^2-2}{6}+\frac{1}{3\sqrt{v^2+1}},-\frac{v}{3}\Big(1-\frac{1}{\sqrt{v^2+1}}\Big)\right), \\
&{\bn}^{t,n}(u,v)=\bt(u,v),\ {\bs}^{t,n}(u,v)=\bn(u,v).
\end{align*}
}
\end{example}
\begin{example}{\rm(A cuspidal edge)
Let $(\bx,\bn,\bs) : \R^2\to \R^3 \times \Delta$ be 
\begin{align*}
\bx(u,v)&=\left(v\cos u-\sqrt{1+v^2}\cos u, v\sin u+\sqrt{1+v^2}\sin u, u-v\sqrt{1+v^2}\right),\\ 
\bn(u,v)&=\left(\cos u-\frac{v}{\sqrt{1+v^2}}\sin u, \sin u+\frac{v}{\sqrt{1+v^2}}\cos u, \frac{1}{\sqrt{1+v^2}}\right), \\ 
\bs(u,v)&=\frac{1}{\sqrt{1+v^2}}\left(-\sin u, \cos u, -v\right).
\end{align*}
Then $\displaystyle \bt(u,v)=\left(-\cos u-\frac{v}{\sqrt{1+v^2}}\sin u, -\sin u+\frac{v}{\sqrt{1+v^2}}\cos u, \frac{1}{\sqrt{1+v^2}}\right)$ and $(\bx,\bn,\bs)$ is a framed surface with the basic invariants
\begin{align*}
\begin{pmatrix}
{a}_1 & {b}_1\\
{a}_2 & {b}_2
\end{pmatrix}
=
\begin{pmatrix}
0 & 2\sqrt{1+v^2}\\
2v & -2
\end{pmatrix}, \ 
\begin{pmatrix}
{e}_1 & {f}_1 & {g}_1\\
{e}_2 & {f}_2 & {g}_2
\end{pmatrix}
=
\begin{pmatrix}
-\frac{1}{\sqrt{1+v^2}} & \frac{2v}{\sqrt{1+v^2}}  & \frac{1}{\sqrt{1+v^2}}\\
\frac{1}{1+v^2} & 0 & -\frac{1}{1+v^2}
\end{pmatrix}.
\end{align*}
Note that $\bx$ at $(0,0)$ is a cuspidal edge by using the criterion of cuspidal edge in \cite{KRSUY}.
It follows that the curvature $C^F$ of $(\bx, \bn, \bs)$ is given by
$$
J^F(u,v)=4v\sqrt{1+v^2}, \ K^F(u,v)=-\frac{2v}{(1+v^2)^{3/2}}, \ H^F(u,v)=\frac{2v^2}{\sqrt{1+v^2}}.
$$
If we take $\lambda^{s,s}(u,v)=-2(1+v^2)$ and $\theta^{s,s}(u,v)=\pi/2$, then condition (\ref{ss-condition}) is satisfied. 
Therefore, $(\bx,\bn,\bs)$ and $({\bx}^{s,s},{\bn}^{s,s},{\bs}^{s,s})$ are $(\bs,\overline{\bs})$-mates where
\begin{align*}
&{\bx}^{s,s}(u,v)=\left(v\cos u+\sqrt{1+v^2}\cos u, v\sin u-\sqrt{1+v^2}\sin u, u+v\sqrt{1+v^2}\right), \\ 
&{\bn}^{s,s}(u,v)=\bt(u,v),\ {\bs}^{s,s}(u,v)= \bs(u,v).
\end{align*}
Since ${\rm det}(\bb(u,v),\bg(u,v))=0$ for all $(u,v) \in U$, if we take $\lambda^{s,n}(u,v)=-1-v^2$, $\theta^{s,n}(u,v)=\pi/4$ and $(u_0,v_0)=(0,0)$, then we have an involute with respect to $\bs$ at $(0,0)$, $\mathcal{I}^{s}(\bx,\bn,\bs)=({\bx}^{s,n},{\bn}^{s,n},{\bs}^{s,n}),$
\begin{align*}
{\bx}^{s,n}(u,v)&=\left(v\cos u, v\sin u, u\right), \ {\bn}^{s,n}(u,v)=\frac{\left(-\sin u, \cos u, -v\right)}{\sqrt{1+v^2}},\\ {\bs}^{s,n}(u,v)&=\frac{\left(-v\sin u, v\cos u, 1\right)}{\sqrt{1+v^2}}.
\end{align*}
The basic invariants of $\mathcal{I}^{s}(\bx,\bn,\bs)$ are given by
\begin{align*}
&\begin{pmatrix}
{a}_1^{s,n} & {b}_1^{s,n}\\
{a}_2^{s,n} & {b}_2^{s,n}
\end{pmatrix}
=
\begin{pmatrix}
\sqrt{1+v^2} & 0\\
0 & 1
\end{pmatrix},\\
&\begin{pmatrix}
{e}_1^{s,n} & {f}_1^{s,n} & {g}_1^{s,n}\\
{e}_2^{s,n} & {f}_2^{s,n} & {g}_2^{s,n}
\end{pmatrix}
=
\begin{pmatrix}
0 & -\frac{1}{\sqrt{1+v^2}} & -\frac{v}{\sqrt{1+v^2}}\\
-\frac{1}{1+v^2} & 0 & 0
\end{pmatrix}.
\end{align*}
Then $\bx^{s,n}$ is a helicoid surface. If we take $\lambda^{n,s}(u,v)=1+v^2$ and $\theta^{n,s}(u,v)=-\pi/4$, then condition (\ref{CsIs-condition}) is satisfied. Therefore, we have a caustic of the framed surface, $\mathcal{C}^{s}(\mathcal{I}^{s}(\bx,\bn,\bs))=({\bx},{\bn},{\bs})$.
If we take $\lambda^{n,s}(u,v)=-1-v^2$ and $\theta^{n,s}(u,v)=\pi/4$, then condition (\ref{CsIs-condition}) is also satisfied. 
Therefore, a caustic of the involute of the framed surface, 
$$
\mathcal{C}^{s}(\mathcal{I}^{s}(\bx,\bn,\bs))
=({\bx}^{n,s}({\bx}^{s,n}, {\bn}^{s,n}, {\bs}^{s,n}),{\bn}^{n,s}({\bx}^{s,n}, {\bn}^{s,n}, {\bs}^{s,n}),{\bs}^{n,s}({\bx}^{s,n}, {\bn}^{s,n}, {\bs}^{s,n}))
$$
is given by
\begin{align*}
&{\bx}^{n,s}({\bx}^{s,n}, {\bn}^{s,n}, {\bs}^{s,n})(u,v)\\
&=\left(v\cos u+\sqrt{1+v^2}\cos u, v\sin u-\sqrt{1+v^2}\sin u, u+v\sqrt{1+v^2}\right), \\
&{\bn}^{n,s}({\bx}^{s,n}, {\bn}^{s,n}, {\bs}^{s,n})(u,v)=\bt(u,v),\\
&{\bs}^{n,s}({\bx}^{s,n}, {\bn}^{s,n}, {\bs}^{s,n})(u,v)=\bs(u,v).
\end{align*}
Thus, we have $\mathcal{C}^{s}(\mathcal{I}^{s}(\bx,\bn,\bs))=({\bx}^{s,s},{\bn}^{s,s},{\bs}^{s,s})$. It follows that $(\bx,\bn,\bs)$ and $\mathcal{C}^{s}(\mathcal{I}^{s}(\bx,\bn,\bs))$ are $(\bs,\overline{\bs})$-mates. Note that $\bx^{s,s}$ at $(0,0)$ is also a cuspidal edge.
}
\end{example}
\begin{example}{\rm(A cuspidal cross-cap)
Let $(\bx,\bn,\bs) : \R^2\to \R^3 \times \Delta$ be
\begin{align*}
\bx(u,v)=\left(u,v^2,uv^3\right),\ \bn(u,v)=\frac{(-2v^3,-3uv,2)}{\sqrt{4v^6+9u^2v^2+4}},\ \bs(u,v)=\frac{(1,0,v^3)}{\sqrt{1+v^6}}.
\end{align*}
Then $\bt(u,v)=(-3uv^4,2(v^6+1),3uv)/\sqrt{4v^6+9u^2v^2+4}\sqrt{1+v^6}$ and $(\bx,\bn,\bs)$ is a framed surface with the basic invariants
\begin{align*}
&\begin{pmatrix}
{a}_1 & {b}_1\\
{a}_2 & {b}_2
\end{pmatrix}
=
\begin{pmatrix}
\sqrt{1+v^6} & 0\\
\frac{3uv^5}{\sqrt{1+v^6}} & \frac{v\sqrt{4v^6+9u^2v^2+4}}{\sqrt{1+v^6}}
\end{pmatrix}, \\ 
&\begin{pmatrix}
{e}_1 & {f}_1 & {g}_1\\
{e}_2 & {f}_2 & {g}_2
\end{pmatrix}\\
&=
\begin{pmatrix}
0 & -\frac{6v\sqrt{1+v^6}}{4v^6+9u^2v^2+4} & 0\\
-\frac{6v^2\sqrt{1+v^6}}{\sqrt{1+v^6}\sqrt{4v^6+9u^2v^2+4}} & \frac{6u(2v^6-1)}{(4v^6+9u^2v^2+4)\sqrt{1+v^6}} & \frac{9uv^3}{(1+v^6)\sqrt{4v^6+9u^2v^2+4}}
\end{pmatrix}.
\end{align*}
It follows that the curvature $C^F$ of $(\bx, \bn, \bs)$ is given by
\begin{align*}
J^F(u,v)&=v\sqrt{4v^6+9u^2v^2+4}, \ K^F(u,v)=-\frac{36v^3}{(4v^6+9u^2v^2+4)^{3/2}}, \\
H^F(u,v)&=-\frac{3u(5v^6-1)}{4v^6+9u^2v^2+4}.
\end{align*}
If we take 
\begin{align*}
&\lambda^{s,t}(u,v)=-\frac{(4v^6+9u^2v^2+4)\sqrt{1+v^6}}{9uv^2+6\sqrt{1+v^6}}, \  \sin\theta^{s,t}(u,v)=-\frac{1}{\sqrt{1+v^2}}, \\ 
&\cos\theta^{s,t}(u,v)=\frac{v}{\sqrt{1+v^2}},
\end{align*}
then condition (\ref{st-condition2}) is satisfied. 
Therefore, the tangential direction framed surface with respect to $\bs$ of the framed surface $(\bx,\bn,\bs)$, $\mathcal{S}^{t}(\bx,\bn,\bs)=({\bx}^{s,t},{\bn}^{s,t},{\bs}^{s,t})$ is given by
\begin{align*}
{\bx}^{s,t}(u,v)&={\bx}(u,v)+\lambda^{s,t}(u,v)\bs(u,v)\\
&=\left(u-\frac{4v^6+9u^2v^2+4}{9uv^2+6\sqrt{1+v^6}}, v^2, v^3\Big(u-\frac{4v^6+9u^2v^2+4}{9uv^2+6\sqrt{1+v^6}}\Big)\right),\\
{\bn}^{s,t}(u,v)&=\frac{v}{\sqrt{1+v^2}}\bt(u,v)+\frac{1}{\sqrt{1+v^2}}\bn(u,v),\\
{\bs}^{s,t}(u,v)&=-\frac{1}{\sqrt{1+v^2}}\bt(u,v)+\frac{v}{\sqrt{1+v^2}}\bn(u,v).
\end{align*}
Moreover, if we take $\lambda^{t,s}(u,v)=-\lambda^{s,t}(u,v)$ and $\theta^{t,s}(u,v)=-\theta^{s,t}(u,v)$, then ${\mathcal{T}}^{s}\circ{\mathcal{S}}^{t}(\bx,\bn,\bs)=(\bx,\bn,\bs).$ Note that $\bx^{s,t}$ at $(0,0)$ is also a cuspidal cross cap.
}
\end{example}

Nozomi Nakatsuyama, 
\\
Muroran Institute of Technology, Muroran 050-8585, Japan,
\\
E-mail address: 25096009b@muroran-it.ac.jp
\\
\\
Masatomo Takahashi, 
\\
Muroran Institute of Technology, Muroran 050-8585, Japan,
\\
E-mail address: masatomo@muroran-it.ac.jp

\end{document}